\documentclass[11pt]{article}
\usepackage{amsmath,amsfonts,latexsym}
\usepackage{amscd,amssymb,amsmath}
\usepackage[all]{xy}            
\usepackage{caption}   
\usepackage{subcaption}
\usepackage{graphicx}
\usepackage[english]{babel}
\usepackage{indentfirst,enumerate}
\usepackage{authblk}
\usepackage[T1]{fontenc}
\usepackage[cp1250]{inputenc}
\usepackage[breaklinks,colorlinks=true, citecolor=blue, linkcolor = red]{hyperref}
\usepackage{breakcites} 
\usepackage{anysize}
\marginsize{3cm}{2cm}{1cm}{1cm}
\newtheorem{theorem}{Theorem}[section]

\newtheorem{fact}[theorem]{Fact}
\newtheorem{lemma}[theorem]{Lemma}

\newtheorem{corollary}[theorem]{Corollary}
\newtheorem{proposition}[theorem]{Proposition}
\newtheorem{definition}[theorem]{Definition}
\newtheorem{remark}[theorem]{Remark}
\newcommand{\hbx}{\hfill$\Box$}

\newcommand{\ue}{\mathrm{e}}

\newcommand{\br}{\mathbb R}

\newcommand{\bz}{\mathbb Z}

\begin{document}

\title{On the interspike-intervals of periodically-driven integrate-and-fire models}
\date{}
\author[1]{Wac\l aw Marzantowicz\thanks{marzan@amu.edu.pl}}
\author[2]{Justyna Signerska\thanks{j.signerska@impan.pl}}
\affil[1]{Faculty of Mathematics and Computer Sci., \, Adam Mickiewicz University of Pozna\'n, ul. Umultowska 87,\, 61-614 Pozna{\'n},
Poland}
\affil[2]{Institute of Mathematics, Polish Academy of Sciences, ul.\'Sniadeckich 8, 00-956 Warszawa, Poland}

\renewcommand\Authands{ and }

\maketitle
\medskip

\begin{abstract}
We analyze properties of the firing map, which iterations give information about consecutive spikes, for periodically driven linear integrate-and-fire  models. By considering locally integrable (thus in general not continuous) input functions, we generalize some results of other authors. In particular we prove theorems concerning continuous dependence of the firing map on the input in suitable function spaces. Using mathematical study of the displacement sequence of an orientation preserving circle homeomorphism,  we provide also a complete description of the regularity properties of the sequence of interspike-intervals and  behaviour of the interspike-interval distribution. Our results allow to explain some facts concerning this distribution observed numerically by other authors.
These theoretical findings are illustrated by carefully chosen computational examples.
\end{abstract}


\section{Introduction}
The scope of this paper are one-dimensional integrate-and-fire (IF) models
\begin{eqnarray}
  \dot{x}&=& F(t,x), \quad F:\mathbb{R}^2\to\mathbb{R}\label{ogolne1}\\
  \lim_{t\to s^+}x(t) &=& x_r \quad \textrm{if}\  x(s)=x_T, \label{ogolne2}
\end{eqnarray}
The dynamical variable  $x(t)$ evolves according to the differential equation (\ref{ogolne1}) as long as it reaches the threshold-value $x=x_{T}$, say at some time $t_1$. Next it is immediately reset to a resting value $x=x_{\textrm{r}}$ and the system continues again from the new initial condition $(x_{\textrm{r}},t_1)$ until possibly next time $t_2$ when the threshold is reached again, etc..  Hybrid systems of this kind are present in neuroscience, where this threshold-reset behaviour is supposed to mimic \emph{spiking} (generation of action potential) in real neurons. Of course, $x_{\textrm{r}}$ and $x_{T}$ could be arbitrary constant values and, moreover, it is possible to consider varying (i.e. time-dependant) threshold and reset, which allows to introduce to the one-dimensional spiking models some other more biologically realistic phenomena (such as refractory periods and threshold modulation \cite{gedeon}). However, often analysis of models with varying threshold and the reset can be reduced through the appropriate change of variables to studying the case of constant $x_{\textrm{r}}$ and $x_{T}$ (see e.g. \cite{brette1}).

Except for the models of neuron's activity IF systems (and circle mappings induced by them in case of periodic forcing) can also be used in modeling of cardiac rhythms and arrhythmias (\cite{arnold1}), in some engineering applications (e.g. electrical circuits of certain type, see \cite{car-hop}) or as models of many other phenomena, which involve accumulation and discharge processes that  occur on significantly different time scales.

For simplicity set $x_r=0$ and $x_T=1$ and suppose that the equation (\ref{ogolne1}) has the property of existence and uniqueness of the solution for every initial condition $(t_0,x_0)\in\mathbb{R}^2$.
\begin{definition}\label{podstawowadefinicja} The firing map for the system (\ref{ogolne1})-(\ref{ogolne2}) is defined as
\begin{displaymath}
\Phi(t):=\inf\{s>t: \ x(s;t,0)\geq 1\}, \ t\in\mathbb{R},
\end{displaymath}
where $x_r=0$, $x_T=1$, and $x(\cdot;t,0)$ denotes the solution of (\ref{ogolne1}) satisfying the initial condition $(t,0)$.
\end{definition}
Of course, the firing map $\Phi(t)$ does not need to be well defined for every $t\in\mathbb{R}$ since for some $t$ it might happen that the solution $x(\cdot;t,0)$ never reaches the value $x=1$. Thus the natural domain of $\Phi$ is the set (compare with \cite{car-ong}):
\begin{displaymath}
\textrm{D}_{\Phi}=\{t\in\mathbb{R}: \ \textrm{there exists}\ s>t  \ \textrm{such that}\ x(s;t,0)=1\}.
\end{displaymath}
Later on we will give necessary and sufficient conditions for the firing map $\Phi:\mathbb{R}\to\mathbb{R}$ of the models considered to be well-defined. \\

The consecutive firing times $t_n$ can be recovered via the iterations of the firing map:
\begin{displaymath}
t_n=\Phi^n(t_0)=\Phi(t_{n-1})=\inf\{s>\Phi^{n-1}(t_0): \ x(s;\Phi^{n-1}(t_0),0)=1\},
\end{displaymath}
and the sequence of interspike-intervals (time intervals between the consecutive resets) is given as
\begin{displaymath}
t_n-t_{n-1}=\Phi^n(t_0)-\Phi^{n-1}(t_0).
\end{displaymath}

There are two basic quantities associated with the integrate-and-fire systems, the firing rate:
\begin{displaymath}
\textrm{FR}(t_0)=\lim_{n\to\infty}\frac{n}{t_n}=\lim_{n\to\infty}\frac{n}{\Phi^n(t_0)},
\end{displaymath}
and its multiplicative inverse, which is the average interspike-interval:
\begin{displaymath}
\textrm{aISI}(t_0)=\lim_{n\to\infty}\frac{t_n}{n}=\lim_{n\to\infty}\frac{\Phi^n(t_0)}{n}.
\end{displaymath}

Obviously, in general the limits above might not exist or depend on the initial condition $(t_0,0)$.

In \cite{keener1} the following observation for periodically driven models was made (the remark was not directly formulated in this way but it is a well-know fact):
\begin{fact}\label{factbykeener} If the function $F$ in (\ref{ogolne1}) is periodic in $t$ (that is, there exists $T$ such that  $F(t,x)=F(t+T,x)$ for all $x$ and $t$), then the firing map $\Phi$ has periodic displacement $\Phi-\textrm{Id}$. In particular for $T=1$ we have $\Phi(t+1)=\Phi(t)+1$ and thus $\Phi$ is a lift of a degree one circle map under the standard projection $\mathfrak{p}: t\mapsto \exp(2\pi \imath t)$.
\end{fact}
In case of periodic forcing, the underlying circle map $\varphi:S^1\to S^1$ such that $\Phi$ is a lift of $\varphi$, is referred as to the \emph{firing phase map}.

Mathematical analysis of one-dimensional IF models was performed e.g. in \cite{brette1,car-ong,gedeon}. Firing map was also investigated combining analytical and numerical approach (\cite{combes} - phase-locking and Arnold tongues, \cite{keener1} - LIF model with sinusoidal
input, \cite{ono} - LIF with periodic input, \cite{tiesinga} -LIF with periodic input and noise, ...). Analytical results concerning the firing map $\Phi$ were obtained assuming that $F(t,x)$ is regular enough  (always at least continuous) and often periodic in $t$.

In particular, we will take into account the Leaky Integrate-and-Fire model (LIF):
\begin{equation}\label{lif}
\dot{x}=-\sigma x +f(t)
\end{equation}  and the Perfect Integrator (PI):
\begin{equation} \label{pi}
\dot{x}=f(t),
\end{equation} where $f:\mathbb{R}\to\mathbb{R}$ will be, usually, periodic and not necessary continuous but only locally-integrable. Allowing also not continuous functions might be important from the point of view of applications where the inputs are often not continuous.  Moreover, although as for the firing map of systems with continuous and periodic drive some rigorous results have been proved (e.g. in \cite{brette1,car-ong,gedeon}), the sequence of interspike-intevals even in such a case, according to our knowledge, has not been investigated in details yet. However, sometimes the sequence of interspike-intevals might be of greater importance than the exact spiking times themselves (\cite{reich}). Interspike-intervals are said to be used  in information encoding by neurons (see \cite{kodowanie} and references therein). Here we will give a detailed description of the sequence of interspike-intervals and interspike-interval distribution with the use of mathematical result concerning  displacement sequence of an orientation preserving circle homeomorphism proved by us in submitted papers ``On the regularity of the displacement sequence of an orientation preserving circle homeomorphism'' and ``Distribution of the displacement sequence of an orientation preserving circle homeomorphism''. However, full exposition of these theorems and proofs in available also in \cite{wmjs4}.
\section{Locally integrable input functions for LIF and PI models: some general properties}

\subsection{Preliminary definitions and facts}
Unless stated otherwise, considering the LIF-model (\ref{lif}) we assume that $\sigma\geq 0$, admitting also $\sigma=0$ to include Perfect Integrator (\ref{pi}) as well. As for the function $f$ in (\ref{lif}) and (\ref{pi}), we assume  that $f\in L^1_{\textrm{loc}}(\mathbb{R})$, i.e. for every compact set $A\subset \mathbb{R}$ the Lebesque integral $\int_{A}\vert f(u)\vert\;du$ exists and is finite. For such functions we redefine the notion of the firing map:
\begin{definition}\label{firingdlalokalnie}
For systems (\ref{lif}) and (\ref{pi}) the firing map $\Phi$ is defined as
\begin{equation}\label{uwiklanenafiringdlalif}
\Phi(t):=\inf\{t_*>t: \ \ue^{\sigma t}\leq \int_t^{t_*}[f(u)-\sigma]\ue^{\sigma u}\;du\}.
\end{equation}
\end{definition}
The above definition is generalization of the ``classical'' firing map $\Phi$ for the differential equation (\ref{lif}) with $f$ continuous, since from Definition \ref{podstawowadefinicja} $\Phi$  has to satisfy the implicit equation:
\begin{equation}\label{implicitnafiringdlalif}
\ue^{\sigma t}=\int_{t}^{\Phi(t)}[f(u)-\sigma]\ue^{\sigma u}\;du.
\end{equation}
\begin{lemma}\label{poprawneokreslenie}
The necessary and sufficient condition for the firing map (\ref{uwiklanenafiringdlalif}) $\Phi:\mathbb{R}\to\mathbb{R}$ to be well-defined is that
\begin{equation}\label{warunek}
\limsup_{t\to\infty}\int_0^t[f(u)-\sigma]\ue^{\sigma u}\;du=\infty
\end{equation}
\end{lemma}
\noindent{\bf Proof.} Suppose that (\ref{warunek}) is satisfied. Choose $t_0\in\mathbb{R}$. Then $\limsup_{t\to\infty}\int_{t_0}^t [f(u)-\sigma]\ue^{\sigma u}\;du=\infty$ and hence there exists $t_*$ such that $\int_{t_0}^{t_*}[f(u)-\sigma]\ue^{\sigma u}\;du\geq \ue^{\sigma t_0}$. Consequently $\Phi(t_0)$ is defined.

Now assume that $\Phi:\mathbb{R}\to\mathbb{R}$ is defined, i.e. for every $t\in\mathbb{R}$ there exists $t_*=\Phi(t)$ such that $\ue^{\sigma t}=\int_{t}^{\Phi(t)}[f(u)-\sigma]\ue^{\sigma u}\;du$. In particular, by Definition \ref{firingdlalokalnie}, taking $t=0$ we obtain that $n=\int_0^{\Phi^n(0)}[f(u)-\sigma]\ue^{\sigma u}\;du$. Thus  $\lim_{n\to\infty}\int_0^{t_n}[f(u)-\sigma]\ue^{\sigma u}\;du=\infty$ for $t_n=\Phi^n(0)$, which proves the statement.   \hbx

\begin{lemma}\label{lematglowny} In the model (\ref{lif}) with $\sigma\geq 0$ and $f\in L^1_{\textrm{loc}}(\mathbb{R})$,  suppose that there exists $\varsigma>0$ such $f(t)-\sigma>\varsigma$ a.e. (i.e. for almost all $t\in\mathbb{R}$ in the sense of Lebesque measure). Then the firing map $\Phi:\mathbb{R}\to \mathbb{R}$ is a homeomorphism.
\end{lemma}
\noindent{\bf Proof.}
Notice that under the stated assumptions, $D_{\Phi}=\mathbb{R}$ on the ground of Lemma \ref{poprawneokreslenie} because for every fixed $t$ the integral $\int_t^{t_*}[f(u)-\sigma]\ue^{\sigma u}\;du$ is a strictly increasing unbounded continuous function of $t_*$. It follows that $\Phi$ is also a continuous monotone function. From  (\ref{uwiklanenafiringdlalif}) we have $0\leq\Phi(t)-t<1/\varsigma$ which gives that $\lim_{t\to\infty}\Phi(t)=\infty$ and $\lim_{t\to - \infty}\Phi(t)=-\infty$ and ends the proof. \hbx

We prove the following simple lemma:
\begin{lemma}\label{skonczonailoscfirings}
Suppose that $f\in L^{1}_{\textrm{loc}}(\mathbb{R})$. Then every run of the model (\ref{lif}) has only finite number of firings in every bounded interval.
\end{lemma}
\noindent{\bf Proof.}
Suppose that there is a firing at time $t_0$. Denote $t_n=\Phi^n(t_0)$ for $n\in\mathbb{N}$. If $\{t_n\}\subset [a,b]$ for some bounded interval $[a,b]\subset \mathbb{R}$ (i.e. $\lim_{n\to\infty}t_n=t_*\in (a,b]$ as sequence $t_n$ is non-decreasing), then from equation (\ref{uwiklanenafiringdlalif}) (or equivalently from the solution $x(t;t_0,0)=\ue^{-\sigma t}\int_{t_0}^{t_1}f(u)\ue^{\sigma u}\;du$ of (\ref{lif}) and the condition $x(t_1;t_0,0)=1$ for the firing at time $t_1$) we obtain that
\begin{displaymath}
\ue^{\sigma t_1}=\int_{t_0}^{t_1}f(u)\ue^{\sigma u}\;du\leq\ue^{\sigma t_1}\int_{t_0}^{t_1}\vert f(u)\vert\;du
\end{displaymath}
and thus $1\leq \int_{t_0}^{t_1}\vert f(u)\vert\;du$ and in general $1\leq \int_{t_{n-1}}^{t_n}\vert f(u)\vert\;du$ for $n\in\mathbb{N}\cup \{0\}$. From this we estimate that
\begin{displaymath}
n\leq \int_{t_0}^{t_n}\vert f(u)\vert\;du\leq \int_{a}^{b}\vert f(u)\vert\;du.
\end{displaymath} As $n$ is arbitrary, it results in $\int_{a}^b\vert f(u)\vert\;du=\infty$, which contradicts that $f\in L^{1}_{\textrm{loc}}(\mathbb{R})$. \hbx
\subsection{Special properties of the Perfect Integrator}
The simple model (\ref{pi}) has many distinct properties than other models. Here we list some of them (for the proofs we refer to \cite{wmjs1}).
\begin{fact}\label{uwagadopisana}
Suppose that $f\in L^{1}_{\textrm{loc}}(\mathbb{R})$ and let $\Phi$ be the firing map  for the Perfect Integrator
(\ref{pi}). Then:
\begin{enumerate}
  \item The consecutive iterates of the firing map are equal to
\begin{equation}\label{szybkieiteracje}
\Phi^n(t)=\min\{s>t: \ x(s;t,0)=n\}
\end{equation}
and  there is only a finite number of firings in every bounded
interval.
  \item $\Phi$ is increasing,
  correspondingly,  non-decreasing, iff
$f(t)>0$, or $f(t)\geq 0$ respectively, a.e. in $\mathbb{R}$.
  \item If $f(t)\geq 0$ a.e., then
\begin{itemize}
\item[i)] {$\Phi$ is left continuous,}
\item[ii)] {$\Phi$  is not right continuous at every point $\bar{a} \in \Phi^{-1}(a)$ for which there exists $\delta_0>0$ such that $f(t)=0$ almost everywhere in $[a,
a+\delta_0]$. Furthermore, such points are the only points of
discontinuity of $\Phi$. }
\end{itemize}
\item If $f(t)>0$ a.e., then $\Phi$ is continuous.
\end{enumerate}
\end{fact}

For the simplified model (\ref{pi}) we even have the analytical expression for the firing rate. Indeed,  the
following theorem was proved in \cite{brette1} (originally for $f$ continuous but the proof is valid for $f\in L^{1}_{\textrm{loc}}(\mathbb{R})$ as well):
\begin{theorem}\label{tw1}
Suppose that for the model  (\ref{pi}) there exists a finite
limit
\begin{equation}\label{row3}
r=\lim_{t\to\infty}\frac{1}{t}\int_{0}^{t}f(u)d u.
\end{equation}
Then for every point $t_0\in\mathbb{R}$ the firing rate $r(t_0)$
exists and is given by the formula (\ref{row3}). In particular, the
firing rate $r(t)$ does not depend on $t$.
\end{theorem}
The proof of the above theorem is immediate: It relies on the fact that  $n=\int^{\Phi^n(t_0)}_{t_0}f(u)\;du$ for every $t_0$ by definition of the firing map and if the above limit exists and equals $r$, then also $\lim_{n\to\infty}\frac{1}{\Phi^n(t_0)}\int^{\Phi^n(t_0)}_{t_0}f(u)\;du=r$.\\
\noindent{\textbf{Example 1.}}
Suppose that the function $f$ in (\ref{pi}) is locally integrable, but is not continuous. Consider, for instance:
\begin{displaymath}
f(t)=\left\{
       \begin{array}{ll}
         2, & \hbox{$t\in[n,n+1/2]$, $\ n\in\bz$;} \\
         0, & \hbox{$t\in(n+1/2,n+1)$.}
       \end{array}
     \right.
\end{displaymath}
In this case we easily get that $\mathcal{M}(f)=1$ (and thus $\varrho=1$). By considering the solution $x(t)$ of $\dot{x}=f(t)$ with the initial condition $x(t_0)=0$, we obtain that
\begin{displaymath}
\Phi(t)=\left\{
          \begin{array}{lll}
            t+1, & \hbox{$t\in(k,k+1/2), \ k\in\bz$;} \\
            k+1/2, & \hbox{$t=k$;} \\
            k+3/2, & \hbox{$t\in[k+1/2,k+1)$}
          \end{array}
        \right.
\end{displaymath}
In particular, $\Phi$ is left-continuous, non-decreasing and constant in the intervals $(k+1/2,k+1)$. However, it is not right-continuous at the points $t=k$. Note that at such points $\Phi(t)=k+1/2$ and $f=0$ in the right neighbourhood $(k+1/2,k+1)$ of $\Phi(t)$ which agrees with Fact  \ref{uwagadopisana} (3.).
\subsection{Continuous dependence on the input function}
\begin{definition}\label{esencjalnesupremum}
The essential supremum of the Lebesque measurable function \newline $f:\mathbb{R}\to\mathbb{R}$ is defined as
\begin{equation}\label{ess}
\textrm{ess sup}f:=\inf\{a\in\mathbb{R}: \ \Lambda(\{t: \ f(t)>a\})=0\},
\end{equation}
where $\Lambda$ denotes the Lebesque measure on $\mathbb{R}$. If $\{a\in\mathbb{R}: \ \Lambda(\{t: \ f(t)>a\})=0\}=\emptyset$, then we write that $\textrm{ess sup} f=\infty$.

If $\textrm{ess sup}\vert f\vert <\infty$, then we say that $f$ is \emph{essentially bounded}.

We also define the essential supremum of $f$ over a compact subset $K\subset\mathbb{R}$ as
\begin{displaymath}
\textrm{ess sup}_K f:=\inf\{a\in\mathbb{R}: \ \Lambda(\{t\in K: \ f(t)>a\})=0\},
\end{displaymath}
\end{definition}
In particular, for the measurable functions $f$ and $g$, $\textrm{ess sup}_K \vert f-g\vert=a_*$ for some $a_*\geq 0$ implies that $\vert f(t)-g(t)\vert\leq a_*$ a.e. in $K$. In general an essentially bounded function does not need to be measurable, since equivalently we might say  that $a_*$ is an essential supremum of $f$ if the set $\{t: \ f(t)>a\}$ is contained in some set of measure zero. However, we will consider only locally integrable functions, thus also measurable. Notice that when $f$ is essentially bounded (and measurable), it is also locally integrable. However, a locally integrable function does not need to be essentially bounded: take for example $f(t)=\frac{1}{\sqrt{\vert t\vert}}$ (with arbitrary finite value at $t=0$).

We consider the space $L^{\infty}_{\textrm{loc}}(\mathbb{R})$ of all locally bounded functions (i.e. $f \in L^{\infty}_{\textrm{loc}}(\mathbb{R})$ iff $\textrm{ess sup}_K \vert f\vert<\infty$ for every compact $K\subset \mathbb{R}$)  as the \emph{Frechet space} with semi-norms and metric defined respectively as
\begin{displaymath}
\|f\|_{L^{\infty}([-k,k])}:=\textrm{ess sup}_{[-k,k]}\vert f\vert
\end{displaymath}
and
\begin{displaymath}
\textrm{d}_{L^{\infty}_{\textrm{loc}}}(f,g):=\sum_{k=1}^{\infty}\frac{1}{2^k}\frac{\|f-g\|_{L^{\infty}([-k,k])}}{1+\|f-g\|_{L^{\infty}([-k,k])}}.
\end{displaymath}
We will mainly consider measurable functions $f\in L^{\infty}_{\textrm{loc}}(\mathbb{R})$. Note that such functions form a subspace of $L^1_{\textrm{loc}}(\mathbb{R})$, which is again a Frechet space with the following semi-norms:
\begin{displaymath}
\|f\|_{L^1([-k,k])}:=\int_{[-k,k]}\vert f(u)\vert\;du, \qquad k=1,2,3,...
\end{displaymath}
The metric on $L^1_{loc}(\mathbb{R})$ can be defined as
\begin{displaymath}
\textrm{d}_{L^1_{\textrm{loc}}}(f,g):=\sum_{k=1}^{\infty}\frac{1}{2^k}\frac{\|f-g\|_{L^1([-k,k])}}{1+\|f-g\|_{L^1([-k,k])}}.
\end{displaymath}
$L^1_{\textrm{loc}}(\mathbb{R})$ with this metric is a complete metric space (see for instance \cite{norma1} p. 2).

Similarly in spaces $C^0(\mathbb{R})$ and $C^m(\mathbb{R})$ of, respectively, continuous and $m$-times continuously differentiable functions $f:\mathbb{R}\to\mathbb{R}$, we introduce the metrics $d_{C^0(\mathbb{R})}$ and $d_{C^m(\mathbb{R})}$ with the use of semi-norms:
\begin{displaymath}
\|f\|_{C^0([-k,k])}:=\textrm{sup}_{t\in [-k,k]}\vert f(t)\vert, \quad k=1,2,3,...
\end{displaymath}
and
\begin{displaymath}
\|f\|_{C^m([-k,k])}:=\max\{\textrm{sup}_{t\in [-k,k]}\vert f(t)\vert,\textrm{sup}_{t\in [-k,k]}\vert f^{(1)}(t)\vert, ..., \textrm{sup}_{t\in [-k,k]}\vert f^{(m)}(t)\vert\},  \quad k=1,2,3,...
\end{displaymath}
where $f^{(n)}(t)$ is the $n$-th derivative of $f$.

\begin{proposition}\label{pierwszee}
In the model (\ref{lif}) with $\sigma\geq 0$ and measurable $f\in L^{\infty}_{\textrm{loc}}(\mathbb{R})$, the mapping $f\mapsto \Phi$ is continuous from the $L^{\infty}_{\textrm{loc}}(\mathbb{R})$-topology into $\mathcal{C}^0(\mathbb{R})$-topology at every point $f$ satisfying $f(t)-\sigma>\varsigma$ a.e. for some $\varsigma>0$.
\end{proposition}
The above proposition says, in particular, that if we have a family of systems $\dot{x}=-\sigma x +f_{\omega}(t)$, where $\omega \in \Omega\subset \mathbb{R}^k$ parameterizes $\{f_{\omega}\}$ continuously in the $L^{\infty}_{\textrm{loc}}(\mathbb{R})$-topology and  $\inf\{\varsigma(\omega): \ f_{\omega}(t)-\sigma>\varsigma(\omega) \ \textrm{a.e.}\}>0$, then the enough small  change of parameter $\omega$ causes an arbitrary small change of the firing map $\Phi$ in the $C^0(\mathbb{R})$-topology (but, of course, even if the firing maps $\Phi_{\omega_1}(t)$ and $\Phi_{\omega_2}(t)$ are uniformly close, the firing times $t^{(1)}_n=\Phi_{\omega_1}^n(t_0)$ and $t^{(2)}_n=\Phi_{\omega_2}^n(t_0)$ with $n\to\infty$ might deviate a lot from each other).

\noindent{\bf Proof of Proposition \ref{pierwszee}.}
Let $f(t)-\sigma>\varsigma>0$ a.e. Our aim is to prove
\begin{equation}\label{najwazniejsze}
\forall_{\varepsilon>0}\ \exists_{\delta>0} \ \forall_{g:\mathbb{R}\to\mathbb{R}} \ d_{L^{\infty}_{\textrm{loc}}(\mathbb{R})}(f,g)<\delta \ \implies\ d_{C^0(\mathbb{R})}(\Phi_f,\Phi_g)<\varepsilon,
\end{equation}
where $\Phi_f$ and $\Phi_g$ are the firing maps induced by $\dot{x}=-\sigma x+f(t)$ and $\dot{x}=-\sigma x +g(t)$, respectively ($f$ and $g$ satisfy requirements stated in Proposition \ref{pierwszee}). Firstly we prove that:
\footnotesize
\begin{equation}\label{posrednie}
\forall_{\varepsilon>0}\ \forall_{K\in\mathbb{N}}\ \exists_{N\in\mathbb{N}}\ \exists_{\delta>0}\ \forall_{g:\mathbb{R}\to\mathbb{R}}\ \underset{[-N,N]}{\textrm{ess sup}}\vert f(t)-g(t)\vert<\delta\ \implies \ \underset{t\in[-K,K]}{\sup}\vert\Phi_f(t)-\Phi_g(t)\vert<\frac{\varepsilon}{2}
\end{equation}
\normalsize
Let then $K\in\mathbb{N}$ be arbitrary. Define $N:=K+\lceil \frac{2}{\varsigma}\rceil$, where $\lceil \frac{2}{\varsigma}\rceil$ is the smallest integer greater or equal $2/\varsigma$. Let $\delta:=\min\{\frac{\varsigma}{2},\frac{\varsigma^2\varepsilon}{4}\}$. Choose the function $g$ satisfying the assumptions and such that $\textrm{ess sup}_{[-N,N]}\vert f(t)-g(t)\vert<\delta$. Then $g(u)-\sigma>\varsigma/2>0$ a.e. in $[-N,N]$. Let then $t\in[-K,K]$ be fixed and suppose that $\Phi_g(t)>\Phi_f(t)$. By definition of the firing map,
\begin{displaymath}
\ue^{\sigma t}=\int_t^{\Phi_f(t)}[f(u)-\sigma]\ue^{\sigma u}\;du=\int_t^{\Phi_f(t)}[g(u)-\sigma]\ue^{\sigma u}\;du.
\end{displaymath}
It follows that
\begin{displaymath}
\int_{\Phi_f(t)}^{\Phi_g(t)}[g(u)-\sigma]\ue^{\sigma u}\;du=\int_t^{\Phi_f(t)}[f(u)-g(u)]\ue^{\sigma u}\;du.
\end{displaymath}
Since $0\leq \Phi_f(t)-t<1/\varsigma$ by the assumption on $f$ and $t, t+1/\varsigma\in [-N,N]$ by our choice of $N$, we estimate
\begin{displaymath}
\int_{\Phi_f(t)}^{\Phi_g(t)}[g(u)-\sigma]\ue^{\sigma u}\;du=\int_t^{\Phi_f(t)}[f(u)-g(u)]\ue^{\sigma u}\;du<\delta(\Phi_f(t)-t)\ue^{\sigma \Phi_f(t)}<\frac{\delta}{\varsigma}\ue^{\sigma \Phi_f(t)}
\end{displaymath}
Simultaneously
\begin{displaymath}
\int_{\Phi_f(t)}^{\Phi_g(t)}[g(u)-\sigma]\ue^{\sigma u}\;du>\frac{\varsigma}{2}\vert\Phi_g(t)-\Phi_f(t)\vert\ue^{\sigma \Phi_f(t)},
\end{displaymath}
provided that $\Phi_g(t)\leq N$. However, suppose that $\Phi_g(t)>N$. Then $\int_t^N[g(u)-\sigma]\ue^{\sigma u}\;du<\ue^{\sigma t}$ by definition of the firing map. On the other hand, by our assumptions on $N$ and $g$,
\begin{displaymath}
\int_t^N[g(u)-\sigma]\ue^{\sigma u}\;du>\frac{\varsigma}{2}(N-t)\ue^{\sigma t}>\frac{\varsigma}{2}\frac{2}{\varsigma}\ue^{\sigma t}=\ue^{\sigma t}
\end{displaymath}  and we arrive at a contradiction. Thus always $\Phi_g(t)\leq N$ and finally we obtain
\begin{displaymath}
\vert \Phi_g(t)-\Phi_f(t)\vert<\frac{2}{\varsigma}\frac{\delta}{\varsigma}\leq \frac{\varepsilon}{2}
\end{displaymath}
If $\Phi_f(t)\geq \Phi_g(t)$, then immediately $\Phi_g(t)\leq N$ (since $\Phi_f(t)\leq N$ as $\Phi_f(t)-t<1/\varsigma$) and we can perform similar calculations.

Now we show how (\ref{posrednie}) implies (\ref{najwazniejsze}). Given $\varepsilon>0$, there exists the smallest integer $K_*$ such that $\sum_{k=K_*}^{\infty}\frac{1}{2^k}\leq\varepsilon/2$ and thus
\begin{displaymath}
\sum_{k=K_*}^{\infty}\frac{1}{2^k}\frac{\|\Phi_f(t)-\Phi_g(t)\|_{C^{0}([-k,k])}}{1+\|\Phi_f(t)-\Phi_g(t)\|_{C^{0}([-k,k])}}<\frac{\varepsilon}{2}
\end{displaymath}
Therefore if also $\sum_{k=1}^{K_*}\frac{1}{2^k}\frac{\|\Phi_f(t)-\Phi_g(t)\|_{C^{0}([-k,k])}}{1+\|\Phi_f(t)-\Phi_g(t)\|_{C^{0}([-k,k])}}<\frac{\varepsilon}{2}$, then $d_{C^{0}(\mathbb{R})}(f,g)<\varepsilon$ (the metric in the Frechet space). But us the function $u\mapsto \frac{u}{1+u}$ is increasing (from $[0,\infty)$ onto $[0,1)$) and the norms $\|\Phi_f(t)-\Phi_g(t)\|_{C^0([-k,k])}$ are non-decreasing with $k$, then

\begin{equation}
\begin{split}
\empty &\sum_{k=1}^{K_*}\frac{1}{2^k}\frac{\|\Phi_f(t)-\Phi_g(t)\|_{C^0([-k,k])}}{(1+\|\Phi_f(t)-\Phi_g(t)\|_{C^0([-k,k])})}\leq \sum_{k=1}^{K_*}\frac{1}{2^k}\frac{\|\Phi_f(t)-\Phi_g(t)\|_{C^0([-K_*,K_*])}}{(1+\|\Phi_f(t)-\Phi_g(t)\|_{C^0([-K_*,K_*])})}\leq \nonumber\\
\empty &\sum_{k=1}^{K_*}\frac{1}{2^k}\|\Phi_f(t)-\Phi_g(t)\|_{C^0([-K_*,K_*])}<\|\Phi_f(t)-\Phi_g(t)\|_{C^0([-K_*,K_*])}.\nonumber
\end{split}
\end{equation}
Now from (\ref{posrednie}) we know that there exists $N_*=K_*+\lceil\frac{2}{\varsigma}\rceil$ and $\widetilde{\delta}$ such that $\|f-g\|_{L^{\infty}([-N_*,N_*])}<\widetilde{\delta}$ implies $\|\Phi_f(t)-\Phi_g(t)\|_{C^0([-K_*,K_*])}<\varepsilon/2$ and thus it also implies $d_{C^{0}(\mathbb{R})}(f,g)<\varepsilon$. But then
\begin{displaymath}
\|f-g\|_{L^{\infty}_{\textrm{loc}}(\mathbb{R})}<\frac{1}{2^{N_*}}\frac{\|f-g\|_{L^{\infty}([-N_*,N_*])}}{1+\|f-g\|_{L^{\infty}([-N_*,N_*])}}<\frac{1}{2^{N_*}}\frac{\widetilde{\delta}}{1+\widetilde{\delta}}
\end{displaymath}
Therefore with $\delta:=\frac{1}{2^{N_*}}\frac{\widetilde{\delta}}{1+\widetilde{\delta}}$ we have $d_{L^{\infty}_{\textrm{loc}}(\mathbb{R})}(f,g)<\delta \ \implies \ d_{C^0_{\textrm{loc}}(\mathbb{R})}(\Phi_f,\Phi_g)<\varepsilon$, which proves the statement. \hbx

Under stronger assumptions on $f$ we  prove the following:
\begin{proposition}\label{ciagloscc1}
If $f\in \mathcal{C}^0(\mathbb{R})$, then  the mapping $f\mapsto \Phi$ is continuous from the topology $\mathcal{C}^0(\mathbb{R})$ into  $\mathcal{C}^1(\mathbb{R})$-topology at every point $f$ such that  there exist $\varsigma>0$ and $M$ with $\varsigma<f(t)-\sigma<M$ for all $t$.
\end{proposition}
\noindent{\bf Proof.}
The equation (\ref{implicitnafiringdlalif}), equivalent to $e^{\sigma t}=H(\Phi(t),t)$ with $H(x,t)=\int_t^x[f(u)-\sigma]\ue^{\sigma u}\;du$, differentiated with respect to $t$ yields that
\begin{equation}\label{pochodnadlalif}
\Phi^{\prime}(t)=\frac{f(t)}{f(\Phi(t))-\sigma}\ue^{-\sigma(\Phi(t)-t)}.
\end{equation}
Note that this formula is well-defined for all $t$ since by our assumption $f(\Phi(t))-\sigma\neq 0$.

Suppose now that $\|f-g\|_{L^{\infty}([-N_*,N_*])}<\delta$ (notation as in the previous proof). Then for $t\in [-K_*,K_*]$ we have the following estimates: $\Phi_f(t)-t<1/\varsigma$ and $\ue^{-\sigma(\Phi_f(t)-t)}<M/\varsigma$, correspondingly $\Phi_g(t)-t<2/\varsigma$  and $\ue^{-\sigma(\Phi_g(t)-t)}<4M/\varsigma$, which can be obtained from (\ref{implicitnafiringdlalif}). Thus
\footnotesize{
\begin{eqnarray}
  \vert \Phi_f^{\prime}(t)-\Phi_g^{\prime}(t) \vert&=& \vert \frac{f(t)}{f(\Phi_f(t))-\sigma}\ue^{-\sigma(\Phi_f(t)-t)}-\frac{g(t)}{g(\Phi_g(t))-\sigma}\ue^{-\sigma(\Phi_g(t)-t)}\vert\leq \nonumber\\
  \empty &\leq& \vert \frac{(f(t)-g(t))(g(\Phi_g(t))-\sigma)\ue^{-\sigma(\Phi_f(t)-t)}}{(f(\Phi_f(t))-\sigma)(g(\Phi_g(t))-\sigma)}\vert + \nonumber\\
  \empty &\empty& \vert \frac{g(t)[\ue^{-\sigma(\Phi_f(t)-t)}(g(\Phi_g(t))-\sigma)-\ue^{-\sigma(\Phi_g(t)-t)}(f(\Phi_f(t))-\sigma)]}{(f(\Phi_f(t))-\sigma)(g(\Phi_g(t))-\sigma)}\vert< \nonumber\\
  \empty &<& \frac{M\delta}{\varsigma^2}+\frac{g(t)[\vert\ue^{-\sigma(\Phi_f(t)-t)}(g(\Phi_g(t))-f(\Phi_f(t)))\vert+\vert (f(\Phi_f(t))-\sigma)(\ue^{-\sigma(\Phi_f(t)-t)}-\ue^{-\sigma (\Phi_g(t)-t)}) \vert]}{\vert(f(\Phi_f(t))-\sigma)(g(\Phi_g(t))-\sigma)\vert} \nonumber\\
  \empty &<& \frac{M\delta}{\varsigma^2} + \frac{2(2M+\sigma)}{\varsigma^2}(\frac{M\delta}{\varsigma}+\frac{4M^2\sigma}{\varsigma}\vert \Phi_f(t)-\Phi_g(t)\vert)\nonumber
\end{eqnarray}}
\normalsize
As $\delta\to 0$ also $\vert\Phi_f(t)-\Phi_g(t)\vert\to 0$ uniformly in $t\in [-K_*,K_*]$ by the previous result. This proves the continuity of $f\mapsto \Phi$ from the Frechet space $C^0(\mathbb{R})$ to the Frechet space $C^1(\mathbb{R})$. \hbx
 \begin{remark} For the Perfect Integrator (\ref{pi}) we obtain that
 \begin{equation}\label{pochodnadlapi}
 \Phi^{\prime}(t)=\frac{f(t)}{f(\Phi(t))}
 \end{equation}
 and thus in the same way (easier) we can prove the statement of Proposition \ref{ciagloscc1} for PI-model under the assumption that $f\in \mathcal{C}^{0}(\mathbb{R})$ and $0<\varsigma<f(t)<M$.
 \end{remark}
 \begin{lemma} For the model $\dot{x}=-\sigma x + f(t)$, $\sigma\geq 0$
  \begin{enumerate}[a)]
  \item if $f\in C^k(\mathbb{R})$, where $k\in\mathbb{N}\cup\{0\}$, and $f(t)-\sigma>0$ for all $t$, then $\Phi\in C^{k+1}(\mathbb{R})$,
  \item if $f\in L^{1}_{\textrm{loc}}(\mathbb{R})$ and $f(t)-\sigma>0$ a.e., then $\Phi\in C^0(\mathbb{R})$.
  \end{enumerate}
 \end{lemma}
\noindent{\bf Proof.} The first part is a direct consequence of the formula (\ref{pochodnadlalif}). The second one follows from the properties of the integral of a locally integrable almost everywhere positive function.  \hbx

\section{Periodic drive for LIF and PI models}
\begin{definition}\label{uogolnionaokresowa} We say that a function $f\in L^{1}_{\textrm{loc}}(\mathbb{R})$ is \emph{periodic}, if there exists $T$ such that $f(t+T)=f(t)$ a.e.
\end{definition}
\begin{remark}
Notice that if $f\in \mathcal{C}^0(\mathbb{R})$ is periodic, then the condition $f(t)-\sigma>\varsigma$ for some $\varsigma>0$ reduces to $f(t)-\sigma>0$. In this case $\Phi:\mathbb{R}\to \mathbb{R}$ is a lift of an orientation preserving circle homeomorphism by the Fact \ref{factbykeener}.
\end{remark}
However, for locally integrable periodic functions the requirement $f(t)-\sigma>\varsigma>0$ a.e. is not equivalent to $f(t)-\sigma>0$ a.e.: Take, for example, $\sigma=1$ and $f(t)=1/n+1$ for $t\in [k-\frac{1}{2^{n-1}},k-\frac{1}{2^n})$, $k\in\mathbb{Z}$, $n\in\mathbb{N}$.  Nevertheless, for $f\in L^1_{\textrm{loc}}(\mathbb{R})$ periodic it is enough to assume that $f(t)-\sigma>0$ a.e. in order to assure that the firing map $\Phi$ (in the generalized sense of Definition \ref{firingdlalokalnie}) has the desired property:
\begin{lemma}\label{uogolnionykeener}
If $f\in L^1_{\textrm{loc}}(\mathbb{R})$ is periodic with period $T=1$ and $f(t)-\sigma>0$ a.e., then the firing map $\Phi$ induced by (\ref{lif}) is a lift of an orientation preserving circle homeomorphism.
\end{lemma}
\noindent{\bf Proof.}
From (\ref{implicitnafiringdlalif}) and periodicity of $f$ we have
\begin{displaymath}
\ue^{\sigma(t+1)}=\int_{t+1}^{\Phi(t+1)}[f(u)-\sigma]\ue^{\sigma u}\;du=\int_{t+1}^{\Phi(t+1)}[f(u-1)-\sigma]\ue^{\sigma u}\;du=\ue^{\sigma}\int_t^{\Phi(t+1)-1}[f(u)-\sigma]\ue^{\sigma u}\;du
\end{displaymath}
which is equivalent to
\begin{displaymath}
\int_{t}^{\Phi(t)}[f(u)-\sigma]\ue^{\sigma u}\;du=\int_t^{\Phi(t+1)-1}[f(u)-\sigma]\ue^{\sigma u}\;du.
\end{displaymath}
Since for fixed $t$, $F(t_*)=\int_t^{t_*}[f(u)-\sigma]\ue^{\sigma u}\;du$ is a continuous increasing function of $t_*$ as the integrand is positive a.e., the above implies that $\Phi(t+1)=\Phi(t)+1$ and thus $\Phi$ has the property of a degree one circle map. Then as $\Phi$ is  continuous and increasing, it must be in fact a  lift of an orientation preserving circle homeomorphism. \hbx

Thus when $f\in L^1_{\textrm{loc}}(\mathbb{R})$ is periodic and $f(t)-\sigma>0$ a.e., the unique firing rate $\textrm{FR}(t)=r$ always exists (independently of $t$), since it is the reciprocal of the unique rotation number $\varrho(\Phi)=\lim_{n\to\infty}\frac{\Phi^n(t)}{n}$, $t\in\mathbb{R}$.

For the simple model (\ref{pi}), where $f\in L^1_{\textrm{loc}}(\mathbb{R})$, $f(t)>0$ a.e. and $f$ is periodic (with period $1$), $\Phi: \mathbb{R}\to\mathbb{R}$ is a lift of an orientation preserving circle homeomorphism $\varphi$ (which is then the  firing phase map), as follows from the lemma above. Moreover, in \cite{brette1} it was proven that in this case $\varphi: S^1\to S^1$ is always conjugated with the rotation $r_{\varrho}$ by $\varrho$ via the homeomorphism $\gamma$ with a lift $\Gamma$ given as
\begin{equation}\label{wzornasprzezeniedlapi}
\Gamma(t):=\frac{\int_0^t f(u)\;du}{\int_0^1 f(u)\;du}, \quad t\in\mathbb{R}
\end{equation}
Formula (\ref{row3}) for the firing rate when $f$ is periodic with period $T=1$ reduces to $r=\int_{0}^1 f(u)\;du$. Thus
\begin{equation}\label{wzornarodlapi}
\varrho=\frac{1}{\int_0^1f(u)\;du}
\end{equation}
is the analytical expression for the rotation number of $\varphi$. One might check by a short direct calculation that indeed we have $\Gamma(\Phi(t))=\Gamma(t)+\varrho$, where $\Gamma$ is continuous, increasing and satisfies $\Gamma(t+1)=\Gamma(t)+1$, i.e. $\gamma$  conjugates $\varphi$ with $r_{\varrho}$.

Observe that an almost everywhere non-negative function $f: \br \to \br $, where $f\in L^1_{\textrm{loc}}(\br)$,
defines a measure $\mu_f$ on $\br$ for which $f$ is the density (the Radon–Nikodym derivative of $\mu_f$),
i.e.
\begin{equation}\label{measure}
\mu_f(A):= \int_A f(u) du
\end{equation}
where $A$ is any measurable (Borel) subset of $\br$. We have the following result:

\begin{proposition}\label{miara niezmiennicza}
Let $f\in L^1_{\textrm{loc}}(\br)$, $f(t)>0$ a.e. be periodic,
$\Phi: \br \to \br $ be the firing map associated with (\ref{pi}),
and $\mu_f$ the associated with $f$ measure.

Then $\mu_f$ is $\Phi $-invariant, i.e. $\Phi$ preserves the measure
$\mu_f$.
\end{proposition}

\textbf{Proof}. We have to prove that $\mu_f(\Phi^{-1}(A))=
\mu_f(A)$ for every measurable set $A\subset \br$. Since $\mu_f(A)=
{\underset{\mathcal{A}}\inf}\,{\underset{i=1}\sum} \,\int_{a_i}^{b_i}f(u)\,du={\underset{\mathcal{A}}\inf}\,{\underset{i=1}\sum}x(b_i)-x(a_i)$, where $\mathcal{A}$ is a cover of $A$ by open intervals $(a_i,b_i)$,
it is enough to show that
$\mu_f(\Phi^{-1}([a,b]))= \mu_f([a,b])$ for every interval $[a,b]
\subset \br$, with $a<b$. Moreover, since in this case $\Phi$ is a homeomorphism, it is equivalent to show that   $\mu_f([\Phi(a), \Phi(b)]) = \mu_f[a,b]$.  We have
$\mu_f([a, b]) = \int_{a}^b f(u) \,du$ and $\mu_f([\Phi(a), \Phi(b)])
= \int_{\Phi(a)}^{\Phi(b)} f(u)\,du$.  By the definition of the firing map,

$$ \int_a^b f(u)\,du  \;=\; \int_{\Phi(a)}^{\Phi(b)} f(u)\,du,$$
which proves the Proposition. \hbx

Throughout the rest of this section we assume that
\begin{enumerate}[1)]
  \item $f$ is measurable and $f\in L^{\infty}_{\textrm{loc}}(\mathbb{R})$ (thus in particular $f\in L^1_{\textrm{loc}}(\mathbb{R})$)
  \item $f$ is periodic (allowing also the case of $f$ constant) in the sense of Definition \ref{uogolnionaokresowa} (with period $T=1$, without the loss of generality)
  \item $f(t)-\sigma>0$ a.e. in $\mathbb{R}$.
\end{enumerate}
Under these assumptions the firing map $\Phi$ is a lift of a circle homeomorphism $\varphi \sim \Phi \ \mod \ 1$. Note that in this case $\Phi$ satisfies $\Phi(t+1)=\Phi(t)+1$ for every $t\in\mathbb{R}$ and thus the compact convergence in the Frechet space $C^0(\mathbb{R})$ ($C^m(\mathbb{R})$)  is equivalent to the uniform convergence (uniform convergence up to $m$-th derivative) on $\mathbb{R}$ because it is enough to consider $\Phi$ and $\Phi_n$ cut to the interval $[0,1]$ (we say that $\Phi_n$ converges compactly to $\Phi$ if for every $K\subset \mathbb{R}$ compact $\lim_{n\to\infty}\sup_{t\in K}\vert \Phi_n(t)-\Phi(t)\vert=0$). In other words, if we admit only $1$-periodic inputs $f$ and $f_n$, then $\sup_{t\in \mathbb{R}}\vert \Phi_n(t)-\Phi(t)\vert\sup_{t\in [0,1]}\vert \Phi_n(t)-\Phi(t)\vert<\varepsilon$ whenever $d_{L^{\infty}_{\textrm{loc}}(\mathbb{R})}(f_n,f)<\delta$ for sufficiently small $\delta$.

Except for the continuity of the mapping $f\mapsto \Phi$ from the $L^{\infty}_{\textrm{loc}}(\mathbb{R})$-topology into $C^0(\mathbb{R})$, we will also need the continuity $f\mapsto \Gamma$, where $\Gamma: \mathbb{R}\to\mathbb{R}$ is the lift of the map $\gamma: S^1\to S^1$ (semi-)conjugating the firing phase map $\varphi: S^1\to S^1$ with the rotation $r_{\varrho}$, where $\varrho=\varrho(\varphi)\in\mathbb{R}\setminus\mathbb{Q}$ is the rotation number of $\varphi$.
 \begin{lemma}
Suppose that $\varrho(\Phi)\in\mathbb{R}\setminus \mathbb{Q}$, where $\Phi$ is a firing map induced by the equation (\ref{lif}) with $\sigma\geq 0$. Then the mapping $f \mapsto \Gamma$, where $\Gamma: \mathbb{R}\to\mathbb{R}$ is a lift of $\gamma$ (semi-)conjugating $\varphi$ with the rotation $r_{\varrho}$, is continuous from the $L^{\infty}_{\textrm{loc}}(\mathbb{R})$-topology into $C^0(\mathbb{R})$ (with $\sup_{\mathbb{R}}$) at every point $f$ such that $f(t)-\sigma>\varsigma$ a.e. for some $\varsigma>0$.
 \end{lemma}
 By the continuity of $f\mapsto \Gamma$ we mean that when $\widetilde{f}$ is a small enough perturbation of $f$, with respect to $L^{\infty}_{\textrm{loc}}$-topology, and $\widetilde{\varrho}=\varrho(\widetilde{\Phi})\in\mathbb{R}\setminus\mathbb{Q}$, then $\widetilde{\Gamma}$ can be chosen such that $\Gamma$ and $\widetilde{\Gamma}$ are uniformly close, where  $\widetilde{\Phi}$ is a firing map induced by $\dot{x}=-\sigma x +\widetilde{f}(t)$ and  $\widetilde{\Gamma}$ is a lift of $\widetilde{\gamma}$ where $\widetilde{\gamma}\circ \widetilde{\varphi}=r_{\widetilde{\varrho}}\circ \widetilde{\gamma}$.

 \noindent{\textbf{Proof.}} We have already proved the continuity of the mapping $\varphi\mapsto \gamma$ from $C^0(S^1)\to C^0(S^1)$ in Theorem 3.2 in \cite{wmjs4}. From this follows the continuity of $\Phi\to\Gamma$ from $C^0(\mathbb{R})$ into $C^0(\mathbb{R})$ (with $\sup_{\mathbb{R}}$-topologies). Since we also have the continuity of $f\to\Phi$ under the stated assumptions, the statement of the above lemma holds.\hbx

\subsection{Regularity properties of the $\textrm{ISI}_n$ sequence}
We will formulate some detailed results concerning regularity of the sequence of interspike-intervals for PI and LIF models. By regularity properties we mean periodicity, asymptotic periodicity and the property of almost strong recurrence.

Due to Lemma \ref{uogolnionykeener}  investigation  of interspike-intervals $\textrm{ISI}_n(t_0)$  for $f$ periodic is covered by the analysis of the displacement sequence $\eta_n(z_0)$ of an orientation preserving circle homeomorphism, being the firing phase map $\varphi$. Thus $\textrm{ISI}_n(t_0)$ equals $\eta_n(z_0)$ (where $z_0=\ue^{2\pi\imath t_0}$) up to some constant integer and the sequences $\textrm{ISI}_n(t_0)$ and $\eta_n(z_0)$ have virtually the same properties.

\begin{proposition}\label{pierwsze} Consider the Perfect Integrator model $\dot{x}=f(t)$. If $T=\int_0^1 f(u)\;du=q/p\in\mathbb{Q}$, then the sequence $\textrm{ISI}_n(t)$ for every initial condition $(t,0)$ is periodic with period $q$.
\end{proposition}
\noindent{\bf Proof.} The analytical expression for the firing rate (\ref{row3}) of PI model for $f$ 1-periodic reduces to $\int_0^1 f(u)\;du$. This means that in our case the rotation number of the underlying firing phase map $\varphi:S^1\to S^1$ equals to $\varrho=1/\int_0^1 f(u)\;du$. Thus if $T=1/\varrho=q/p$ is rational, $\varphi$ is topologically conjugated to the rational rotation $r_{\varrho}$ by $\varrho$  and thus there are only periodic orbits with period $q$. As a  result, the sequence of displacements of $\varphi$, and consequently the sequence $\textrm{ISI}_n(t_0)$, is periodic with period $q$. \hbx

\noindent{\textbf{Example 2.}} For the LIF model $\dot{x}=-x+\frac{1}{1-\ue^{-q}}$, where $q\in\mathbb{N}$, the sequence of interspike-intervals is constant: $\textrm{ISI}_n(t_0)=q$. Indeed, in \cite{brette1} it was shown that the input current of such a form induces conjugacy with the rational rotation by $\varrho(\Phi)=q$. Consequently, the firing map $\Phi$ satisfies $\Phi(t)=t+q$ and is simply a translation by $q$. Thus for every $n\in\mathbb{N}$ and $t$ we have $\Phi^n(t)-\Phi^{n-1}(t)=q$ and we observe 1 spike per every $q$ periods of forcing.

R. Brette (\cite{brette1}) also proved that $f(t)=\frac{1}{1-\ue^{-q}}$ is the only one input current which induces conjugacy with a rotation by $q\in\mathbb{N}$ (for $\sigma=1$). It is much harder to show what are all the input currents that induce conjugacy with $\varrho=q/p$ ($p\neq 1$) but this assumption implies some constraints on $f(t)$, which seem to be quite restrictive (for some values of $p/q$ the conjugacy might not be possible at all, see discussion in \cite{brette1}). Thus we might conclude that in ``majority of cases'' the firing phase map arising from the LIF model, which has rational firing rate, is not conjugated to the corresponding rotation and:
\begin{remark}\label{drugie} For the LIF model with a firing rate $\textrm{FR}=q/p$, the sequence of interspike-intervals $\textrm{ISI}_n(t_0)$ is ``typically'' not periodic but only asymptotically periodic (with the period equal to $q$ in the limit $n\to\infty$). Precisely,
\begin{displaymath}
\forall_{\varepsilon>0} \ \exists_{N\in\mathbb{N}} \ \forall_{n\in\mathbb{N}} \ \forall_{k\in\mathbb{N}} \ \vert \textrm{ISI}_{n+kq}(t_0)-\textrm{ISI}_{n}(t_0)\vert<\varepsilon
\end{displaymath}
\end{remark}
\noindent{\bf Proof.} This is a direct consequence of Proposition 2.5 in \cite{wmjs4}. \hbx

\begin{figure}[h!]
\centering
                \includegraphics[width=0.5\textwidth]{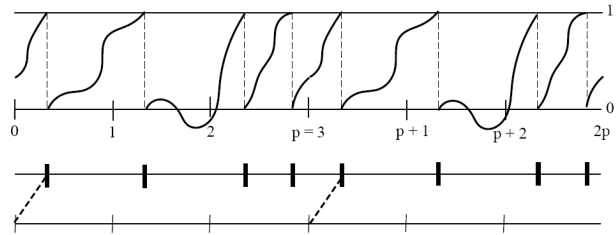}
                \caption[An example of a phase-locking state.]{An example of $4:3$-phase locking}\label{fig:locking}
\end{figure}

When the input function is periodic, the (asymptotically) periodic output of the system, in terms of interspike-intervals, is connected with the phenomena called phase-locking, see Figure \ref{fig:locking}. Precisely, we say that the system exhibits $q:p$ - phase locking (which corresponds to the rotation number equal to $p/q$), when it fires $q$ spikes for every $p$ cycles of forcing (the spikes occur in fixed phases of the forcing period) and this state is structurally stable, i.e. it persists under a small change of a parameter $\theta\in\Theta$. Types of phase-locking change with the change of the rotation number, but the mapping $\theta\mapsto\varrho$ is usually constant (under some conditions) on rational values of $\varrho$ (look for such concepts as the devil-staircase and Arnold-tongues).

In next we pass to the case of irrational firing rate. The same property as for the displacement sequence of a circle homeomorphism with the irrational rotation number can be shown for the sequence of interspike-intervals for the LIF model (see Proposition 3.20 in \cite{wmjs4}):

\begin{theorem}\label{isi strongly recurrent} Consider the LIF model $\dot{x}=-\sigma x +f(t)$ ($\sigma\geq 0$) where $f\in L^1_{\textrm{loc}}(\mathbb{R})$ is periodic, $f(t)-\sigma>0$ a.e. and the rotation number $\varrho(\Phi)$ is irrational.

Then the sequence $\{\textrm{ISI}_n(t_0)\}$ is almost strongly recurrent for all $t_0\in\widetilde{\Delta}$, where $\widetilde{\Delta}$ is a lift to $\mathbb{R}$ of the underlying minimal set $\Delta\subset S^1$ (possibly $\Delta = S^1$), i.e.
\begin{displaymath} \forall_{\varepsilon>0} \ \exists_{N\in\mathbb{N}} \ \forall_{n\in N} \ \forall_{k\in\mathbb{N} \cup \{0\}} \ \exists_{i\in\{0,1,..., N\}} \ |\textrm{ISI}_{n+k+i}(t_0)-\textrm{ISI}_{n}(t_0)|<\varepsilon
\end{displaymath}
Moreover, if $f\in C^2(\mathbb{R})$, then the sequence $\{\textrm{ISI}_n(t_0)\}$ is almost strongly recurrent for all $t_0\in\mathbb{R}$ (in this case $\Delta =S^1$).
\end{theorem}
\noindent{\bf Proof.} Under the stated assumptions $\varphi:S^1\to S^1$ is a homeomorphism with irrational rotation number. For $t_0\in\widetilde{\Delta}$ the underlying displacement sequence $\eta_n(z_0)=\Phi^n(t_0)-\Phi^{n-1}(t_0) \ \mod \ 1$, $z_0=\ue^{2\pi\imath t_0}$, is almost strongly recurrent by Proposition 3.20 in \cite{wmjs4}. But then the sequence of interest $\textrm{ISI}_n(t_0)=\Phi^n(t_0)-\Phi^{n-1}(t_0)$ is almost strongly recurrent as well.

As for the second part of the statement, notice that if $f\in C^2(\mathbb{R})$ , then $\varphi\in C^2(S^1)$ and thus on the ground of the Denjoy Theorem (\cite{denjoy}), $\varphi$ is transitive and $\{\textrm{ISI}_n(t_0)\}$ is almost strongly recurrent for all $t_0\in\mathbb{R}$. \hbx\\

\begin{remark} Notice that Proposition \ref{pierwsze}, Remark \ref{drugie} and Theorem \ref{isi strongly recurrent}  remain true when one replaces the ``sequence of interspike-intervals'' simply with the ``sequence of firing times phases'', i.e. the sequence $\Phi^n(t_0) \ \mod \ 1$.
\end{remark}
\subsection{Distribution of interspike-intervals}
In this part we consider IF models, for which the firing rate, and consequently the rotation number of the firing phase map $\varphi$, is irrational.

\begin{proposition}\label{gestosc w przedziale dla lif}
Consider the LIF model $\dot{x}=-\sigma x +f(t)$, where $f\in L^1_{\textrm{loc}}(\mathbb{R})$ is $1$-periodic, $f(t)-\sigma>0$ a.e. and the rotation number $\varrho=\varrho(\Phi)$ is irrational. By $\Gamma$ denote the lift of $\gamma$ (semi-)conjugating $\varphi$ with $r_{\varrho}$. Under these assumptions
\begin{enumerate}[1)]
  \item If $\varphi$ is transitive (for example when $f\in C^2(\mathbb{R})$), then the sequence $\textrm{ISI}_n(t_0)$ for every $t_0\in\mathbb{R}$ is dense in the interval
  \begin{equation} \label{wzor na przedzial}
  \mathcal{S}=\Psi([0,1])=\Omega([0,1]),
  \end{equation}
  where $\Psi(t)=\Phi(t)-t$ is a displacement function of $\Phi$ and $\Omega(t):=\Gamma^{-1}(t+\varrho)-\Gamma^{-1}(t)$.
  \item If $\varphi$ is not transitive, then the sequence $\textrm{ISI}_n(t_0)$ for $t_0\in\widetilde{\Delta}$ (the total lift of $\Delta$ to $\mathbb{R}$) is dense in the set
   \begin{displaymath}
  \widehat{\mathcal{S}}=\Psi(\widetilde{\Delta})=\widehat{\Omega}(\widehat{\Delta}_0),
  \end{displaymath}
  where  $\widehat{\Delta}_0$ is a lift to $[0,1]$ of a subset $\widehat{\Delta}\subset \Delta$, such that the semi-conjugacy $\gamma$ is invertible on $\widehat{\Delta}$, and $\widehat{\Omega}:=\widehat{\Gamma}^{-1}(t+\varrho)-\widehat{\Gamma}^{-1}(t)$, where $\widehat{\Gamma}:=\Gamma\upharpoonright\widehat{\Delta}_0$ is a lift of $\gamma$ cut to the set $\widehat{\Delta}_0$. \\
  Moreover, when one takes $t_0\in\mathbb{R}\setminus\widetilde{\Delta}$, then for every $t\in\widetilde{\Delta}$ there exists an increasing sequence $n_k$, $k\in\mathbb{N}$, such that for every $l\in\mathbb{N}$ we have
  \begin{displaymath}
  \lim_{k\to\infty}\textrm{ISI}_l(\Phi^{n_k}(t_0))=\Phi^{n_k+l}(t_0)-\Phi^{n_k+l-1}(t_0) =\Phi^l(t)-\Phi^{l-1}(t)=\textrm{ISI}_l(t)
  \end{displaymath}
\end{enumerate}
\end{proposition}

Since usually we do not know the formula for the (semi-)conjugacy $\Gamma$ (except for the Perfect Integrator), the equivalent formula for the concentration set of $\textrm{ISI}$ involving $\Omega$ is not directly useful but it is used in proving statements concerning the  distribution of interspike-intervals with respect to the unique invariant measure.

\noindent{\bf Proof of Proposition \ref{gestosc w przedziale dla lif}.} The above proposition is a direct consequence of Proposition 3.1 in \cite{wmjs4}, the corresponding statement for the displacement sequence of an orientation preserving circle homeomorphism with irrational rotation number.  \hbx

 Proposition \ref{gestosc w przedziale dla lif} means that even if we are not able to compute directly the set of concentration of interspike-intervals, we know at least that interspike-intervals practically fill a whole interval (i.e. do not form for instance something in the type of a Cantor set), if only $f$ is smooth enough. This is also visible in numerical Example 5. However, note that in a special case, where $\varphi$ is a strict rotation (as happens for example for LIF and PI with constant input), this interval degenerates to a single point $\{\varrho\}$ (since the rotation is an isometry).

We will discuss the distribution $\mu_{\textrm{ISI}}$ of interspike-intervals with respect to the unique (up to normalization), $\varphi$-invariant measure $\mu$.
\begin{definition} Suppose that the rotation number $\varrho(\Phi)$ is irrational. Let $\mu$ be the unique invariant probability measure for $\varphi\sim\Phi\mod 1$. The distribution of interspike-intervals is defined as
\begin{displaymath}
\mu_{\textrm{ISI}}(A):=\mu(\{t\in [0,1]:\ \Phi(t)-t\in A\})=\mu(\Psi^{-1}(A)), \quad A\subset\mathbb{R}
\end{displaymath}
where $\Psi(t)=\Phi(t)-t$, $t\in [0,1]$, is a displacement function associated with $\Phi$.
\end{definition}
Note that since $\Phi \ \mod \ 1$ is periodic with period $1$, we consider only $t\in[0,1]$. Moreover, although the measure $\mu$ has support contained in $[0,1]$, as it is the invariant measure for $\Phi \mod 1: [0,1]\to [0,1]$, the measure $\mu_{\textrm{ISI}}$ has support equal to $\Psi([0,1])$, which in general might not be contained in $[0,1]$ but it is always contained is some interval of length not greater than $1$ because $\Phi$ maps intervals of length $1$ onto intervals of length $1$ due to the fact that $\Phi(t+1)=\Phi(t)+1$ for every $t$ (the resulting interval, containing $\textrm{supp}(\mu_{ISI})$,  is shifted by $a$ from its projection $\mod 1$ into $[0,1]$, where $a>0$ is such that $\Phi(0)=a$). We can consider $\mu_{\Psi}(A)$, where $A$ is an arbitrary subset of $\mathbb{R}$, if we adopt the convention that $\mu_{\textrm{ISI}}$ is defined on the whole $\mathbb{R}$ but it simply vanishes  everywhere outside ist support.

By $\Lambda$ denote the Lebesque measure on $[0,1]$.
\begin{proposition}\label{postac miary}
Under the assumptions of Proposition \ref{gestosc w przedziale dla lif} the distribution $\mu_{\textrm{ISI}}$ is the transported Lebesque $\Lambda$ measure:
\begin{enumerate}
  \item If $\varphi$ is transitive, then
  \begin{displaymath}
  \mu_{\textrm{ISI}}(A)=\Lambda(\Omega^{-1}(A)),  \ A\subset\mathbb{R}
  \end{displaymath}
  and the support of $\mu_{\textrm{ISI}}$ equals
  \begin{displaymath}
  \textrm{supp}(\mu_{\textrm{ISI}})=\Psi([0,1])=\Omega([0,1]).
  \end{displaymath}
  \item If $\varphi$ is not transitive, then analogously
  \begin{displaymath}
  \mu_{\textrm{ISI}}(A)=\Lambda(\widehat{\Omega}^{-1}(A)), \ A\subset\widehat{S},
  \end{displaymath}
  and the support of $\mu_{\textrm{ISI}}$ equals
  \begin{displaymath}
  \textrm{supp}(\mu_{\textrm{ISI}})=\Psi(\widetilde{\Delta}),
  \end{displaymath}
  where $\Omega$, $\widehat{\Omega}$, $\widetilde{\Delta}$ and $\widehat{S}$ are as in Proposition \ref{gestosc w przedziale dla lif}.
\end{enumerate}
\end{proposition}
Obviously, the support $\textrm{supp}(\mu_{\textrm{ISI}})$ is just the set of concentration of the interspike-intervals sequence.

\noindent{\bf Proof of Proposition \ref{postac miary}.} Notice that  $\mu_{\textrm{ISI}}$ corresponds to the distribution $\mu_{\Psi}$ of the displacement sequence of $\varphi$ and compare with Proposition 3.3 in \cite{wmjs4}.  $\empty\quad \empty$\hbx

We are concerned with the distribution $\mu_{\textrm{ISI}}$ of interspike-intervals with respect to the natural invariant measure $\mu$, since this theoretical distribution is a limiting distribution of interspike-intervals computed along an arbitrary trajectory:
\begin{proposition}\label{srednia}
Under the assumptions of Proposition \ref{gestosc w przedziale dla lif} (regardless the transitivity of $\varphi$), for $A\subset \mathbb{R}$ we have
\begin{displaymath}
\lim_{n\to\infty}\frac{\sharp \{0\leq i \leq n-1: \ \textrm{ISI}_i(t)\in A\}}{n}= \mu_{\textrm{ISI}}(A),
\end{displaymath}
where $\sharp$ denotes the number of elements of the set, and the above convergence is uniform (with respect to $t\in\mathbb{R}$).

The average interspike interval $\textrm{aISI}$ (which equals the rotation number $\varrho(\Phi)$) is the mean of the distribution $\mu_{\textrm{ISI}}$:
\begin{displaymath}
\textrm{aISI}=\int_{\mathbb{R}}\Phi(t)-t\;d\mu(t)=\int_{\mathbb{R}}\;d\mu_{\textrm{ISI}}.
\end{displaymath}
\end{proposition}
\noindent{\bf Proof.} The statements can be justified by the Birkhoff Ergodic Theorem applied to the observable $\Psi(t)=\Phi(t)-t$. The uniform convergence in the first part follows from the fact that the system $(\Phi\mod 1, [0,1], \mu)$ is not only ergodic but uniquely ergodic (compare with Proposition 4.1.13 and Theorem 11.2.9 in \cite{katok}   or Proposition 3.12 in  \cite{wmjs4}).\hbx

Our aim will be to consider parameter-dependant IF systems and to formulate results describing how the distribution $\mu_{\textrm{ISI}}$ varies with change of the parameter.

\begin{definition}[see e.g. \cite{bartoszynski}] Let $X$ be a complete separable metric space and $\mathcal{M}(X)$ the space of all finite measures defined on the Borel $\sigma$-field $\mathcal{B}(X)$ of subsets of $X$.

A sequence $\mu_n$ of elements of $\mathcal{M}(X)$ is called weakly convergent to $\mu\in\mathcal{M}(X)$ if for every bounded and continuous function $f$ on $X$
\begin{displaymath}
\lim_{n\to\infty} \int\limits{X} f(x)\;d\mu_n(x)=\int\limits{X} f(x)\;d\mu(x).
\end{displaymath}
We denote the weak convergence as $\mu_n\implies \mu$.
\end{definition}
\begin{definition} A Borel set $A$ is said to be a continuity set for $\mu$ if $A$ has $\mu$-null boundary, i.e.
\begin{displaymath}
\mu(\partial A)=0
\end{displaymath}
\end{definition}
One can show (cf. \cite{ranga}) that $\mu_n\implies \mu$ if and only if for each continuity set $A$ of $\mu$, $\lim \mu_n(A)=\mu(A)$.
\begin{proposition}\label{slabazb}
Consider the systems $\dot{x}=-\sigma x+f(t)$ and $\dot{x}=-\sigma x+f_n(t)$, $n\in\mathbb{N}$, where the functions $f_n, f \in L^{\infty}_{\textrm{loc}}(\mathbb{R})$ are measurable, periodic with period $1$ and  $f$ satisfies $f(t)-\sigma>\varsigma>0$ a.e.. Suppose that all the induced firing  maps $\Phi_n$ and $\Phi$ have irrational rotation numbers, $\varrho_n$ and $\varrho$, respectively. By $\mu^{(n)}_{\textrm{ISI}}$ and $\mu_{\textrm{ISI}}$ denote the interspike-interval distributions, correspondingly for  $\Phi_n$ and $\Phi$, with respect to the corresponding invariant measures $\mu^{(n)}$ and $\mu$.

If $f_n\to f$ in $L^{\infty}_{\textrm{loc}}(\mathbb{R})$-topology, then
\begin{displaymath}
\mu^{(n)}_{\textrm{ISI}}\implies \mu_{\textrm{ISI}}.
\end{displaymath}
\end{proposition}
\noindent{\bf Proof.} Recall that the invariant measures, $\mu^{(n)}$ and $\mu$, are the Lebesque measure transported by the maps $\Gamma_n$ and $\Gamma$ (semi-)conjugating corresponding firing maps $\Phi_n$ and $\Phi$ with the rotations. Since we already know  that the mapping $f\mapsto \Gamma$ is continuous from the $L^{\infty}_{\textrm{loc}}(\mathbb{R})$-topology into $C^0(\mathbb{R})$, it must hold that $\Gamma_n\to\Gamma$ in $C^0(\mathbb{R})$ (with $\sup_{\mathbb{R}}$-topologies). But then $\mu^{(n)}\implies \mu$, i.e. we have the weak convergence of the invariant measures. Since the interspike-interval distributions are in turn the invariant measures transported by the corresponding displacement functions $\Psi_n\to \Psi$ in $\sup_{\mathbb{R}}$, by the same argument we get the statement on $\mu^{(n)}_{\textrm{ISI}}$ and $\mu_{\textrm{ISI}}$.\hbx\\

 Recall that the weak convergence of measures does not imply the point-wise convergence of the corresponding densities: in general the densities of $\mu^{(n)}_{\textrm{ISI}}$ or $\mu_{\textrm{ISI}}$ might not exist, as in Example 3.

Notice that in the above proposition we needed the fact that all the rotation numbers $\varrho_n$ and $\varrho$ are irrational since this guarantees that the unique invariant measures $\mu^{(n)}$ and $\mu$ exist and we can define the distributions $\mu^{(n)}_{\textrm{ISI}}$ and $\mu_{\textrm{ISI}}$. Later on we will see what happens if the intermediate firing maps $\Phi_n$ may have rational rotation numbers as well.

However, now we want to formulate some more detailed theorems on convergence of interspike-interval distributions. This is quite simply achievable for the simplest model, the Perfect Integrator Model.
\begin{proposition}\label{simplePIM}
Consider the systems of Perfect Integrators $\dot{x}=f_n(t)$, $n\in\mathbb{R}$, and $\dot{x}=f(t)$, where the functions $f_n, f\in C^{0}(\mathbb{R})$ are periodic with period $1$, $f_n(t), \ f(t)>0$, $f_n\to f$ in $C^0(\mathbb{R})$ and where all the rotation numbers of the firing maps are irrational, $\varrho_n,\varrho\in\mathbb{R}\setminus \mathbb{Q}$. Then the invariant measures $\mu^{(n)}$ and $\mu$ have densities, say $g_n$ and $g$ correspondingly,  and
\begin{displaymath}
g_n\to g \quad \textrm{in} \ C^0([0,1]).
\end{displaymath}
As for the distributions of interspike intervals, if additionally the set of critical points of the displacement function $\Psi$ of the limiting firing map $\Phi$, i.e. the set $\mathcal{C}_{\Psi}:=\{t\in [0,1]:\ \Psi^{\prime}(t)=0\iff \Phi^{\prime}(t)=1\}$, is of Lebesque measure $0$, then we have
\begin{equation}\label{jednostajnanaodcinkach}
\sup[\vert\mu^{(n)}_{\textrm{ISI}}(I)-\mu_{\textrm{ISI}}(I)\vert, \ I\in \mathcal{J}]\to 0,
\end{equation}
where $\mathcal{J}$ denotes the class of all intervals $I\subset [0,1]$ (open, closed, half-closed).
\end{proposition}
\noindent{\bf Proof.} Proposition \ref{miara niezmiennicza} provides the following formula for the unique invariant probability measure $\mu$ of the firing phase map for the Perfect Integrator:
\begin{equation}\label{niezmdlapi}
\mu(A)=\frac{\int_{A}f(u)\;du}{\int_{[0,1]}f(u)\;du}, \ A\subset [0,1]\ \textrm{- Borel subset}.
\end{equation}
This formula is also consistent with the formula for the (semi\nolinebreak[10]-\nolinebreak[10])conjugacy $\Gamma$, since by the standard result on circle homeomorphisms (see \cite{one-dim}, p.34) it holds that $\Gamma(t)=\mu([0,t])$ (provided that $\Gamma(0)=0$ which we can assume without the loss of generality as the conjugacy $\Gamma$ is given up to the additive constant). Thus $g_n(t)=f_n(t)/\int_{0}^{1}f_n(u)\;du=\Gamma_n^{\prime}(t)$, $t\in [0,1]$, correspondingly  $g(t)=f(t)/\int_{0}^{1}f(u)\;du=\Gamma^{\prime}(t)$, and the uniform convergence of densities follows. This in particular implies that
\begin{equation}\label{pss}
\sup[\vert\mu^{(n)}(A)-\mu(A)\vert; \ A\textrm{ - Borel subset of }[0,1]]\to 0.
\end{equation}
Indeed, by the existence and convergence of densities $\Gamma_n^{\prime}$  and $\Gamma^{\prime}$, $\vert\mu^{(n)}(A)-\mu(A)\vert=\vert\int_A\Gamma_n^{\prime}(u)-\Gamma^{\prime}(u)\;du\vert<\varepsilon \Lambda(A)$ for sufficiently large $n$, but since $\Lambda(A)\leq 1$, the convergence is uniform with respect to the choice of $A\subset [0,1]$.

Unfortunately, without the assumption on the $0$-measure of the set of critical points of the limiting displacement function $\Psi$, we cannot assure that  the distribution $\mu_{\textrm{ISI}}$ of interspike-intervals has density (with respect to the Lebesque measure), as we will see in Example 3. However, under this assumption by Proposition 3.7 in \cite{wmjs4} we obtain (\ref{jednostajnanaodcinkach}), i.e.  we know that convergence of interspike-intervals distributions is uniform on the collection of all the intervals (in this case the density of $\mu_{\textrm{ISI}}$ exists by Theorem \ref{wzor na gestosc} below, but the densities of $\mu^{(n)}_{\textrm{ISI}}$ might still not exist and we cannot argue as above for (\ref{pss})). \hbx

In order to compute the set $\mathcal{C}_{\Psi}$ for Perfect Integrator one has to solve in $t$ the implicit equation $f(t)=f(\Phi(t))$ by (\ref{pochodnadlapi}), which usually is difficult. But in the forthcoming examples we will see that verifying the assumption on the zero Lebesque measure of this set is sometimes not that challenging. We remark only that  when in particular $f$ is constant, this assumption is not satisfied but in this case the emerging firing map $\Phi$ is exactly the lift of the rotation by $\varrho$ and  the distribution of interspike intervals equals simply the Dirac delta $\delta_{\varrho}$ and thus $\mu_{\textrm{ISI}}$ is not absolutely continuous with respect to $\Lambda$.

 The theorem below provides sufficient conditions for the distribution $\mu_{\textrm{ISI}}$ to have the density with respect to the Lebesque measure. We formulate it  in the most general form:
\begin{theorem}\label{wzor na gestosc} Suppose that the firing  map $\Phi$ arising from the system $\dot{x}=F(t,x)$ is a $C^1$-diffeomorphism with irrational rotation number $\varrho$, which is conjugated with the translation by $\varrho$ via a $C^1$-diffeomorphism $\Gamma$ and that the set $\mathcal{C}_{\Psi}\subset [0,1]$ of critical points of the displacement function $\Psi\upharpoonright_{[0,1]}$ is of Lebesque measure $0$.

Then the distribution $\mu_{\textrm{ISI}}$ is absolutely continuous with respect to the Lebesque measure $\Lambda$ with the density $\Delta(y)$ equal to
\begin{equation}\label{wzor na gestosc_eq}
\Delta(y)=\left\{
            \begin{array}{ll}
              0 & \hbox{if $y\not\in \mathrm{supp}(\mu_{\textrm{ISI}})$;} \\
              \sum_{t\in \Psi^{-1}(y)}\Gamma^{\prime}(t)\frac{1}{\vert\Phi^{\prime}(t)-1\vert} & \hbox{if $y\in \mathrm{supp}(\mu_{\textrm{ISI}})$.}
            \end{array}
          \right.
\end{equation}
where the latter is well-defined almost everywhere in $\mathrm{supp}(\mu_{\textrm{ISI}})$, i.e. in $\mathrm{supp}(\mu_{\textrm{ISI}})\setminus V(C_{\Psi})$, where $V(C_{\Psi})$ denotes the set of critical values of $\Psi\upharpoonright [0,1]$.
\end{theorem}
\noindent{\bf Proof.} The theorem is a mere tautology of Theorem 3.10 in \cite{wmjs4}.\hbx\\

In particular for the Perfect Integrator we make use of the formulas (\ref{pochodnadlapi}) and (\ref{wzornasprzezeniedlapi}) for the derivative $\Phi^{\prime}(t)$ and the conjugacy $\Gamma$ in order to obtain that (\ref{wzor na gestosc_eq}) reduces to:
\begin{equation}
\Delta(y)=\left\{
            \begin{array}{ll}
              0 & \hbox{if $y\not\in [\min_{u\in[0,1)}\Phi(u)-u,\max_{u\in[0,1)}\Phi(u)-u]$;} \nonumber\\
              \sum_{t\in \Psi^{-1}(y)}\frac{f(t)}{\int_0^1f(u)\;du}\frac{f(\Phi(t))}{\vert f(t)-f(\Phi(t))\vert} & \hbox{if $y\in [\min_{u\in[0,1)}\Phi(u)-u,\max_{u\in[0,1)}\Phi(u)-u]$.}\nonumber
            \end{array}
          \right.
\end{equation}

\noindent{\textbf{Example 3.}} Consider the systems $\dot{x}=f_n(t)$, where
\[f_n(t)=A_n+B_n\cos(2\pi n t) \]
 and
\[A_n\to A_0>0 \quad \textrm{and} \quad 0<B_n\to 0.\] Suppose that the constants $A_n$ and $B_n$ are such that $f_n(t)>0$, at least for sufficiently large $n\in\mathbb{N}$. In particular, we have the convergence $f_n\to f_0$ in $C^{1}(\mathbb{R})$, where $f_0\equiv A_0$. The firing maps $\Phi_n$ are then the lifts of circle diffeomorphisms with rotation numbers $\varrho_n=\frac{1}{A_n}$. Moreover, on the account of Proposition \ref{ciagloscc1}, $\Phi_n\to \Phi_0$ in $C^1(\mathbb{R})$, where $\Phi_0(t)=t+\varrho_0=t+\frac{1}{A_0}$ is the firing map induced by the equation $\dot{x}=f_0(t)$ and simply a lift of the rotation by $\varrho_0$. The firing maps $\Phi_n$ and $\Phi_0$ are conjugated to the corresponding rotations, respectively, via diffeomorphisms
\[\Gamma_n(t)=\varrho_n\int_{0}^t f_n(u)\,du=t+\frac{B_n}{2\pi n A_n}\sin(2\pi n t) \]
and
\[ \Gamma_0(t)=t\]
Assume that $\varrho_n, \varrho_0\in\mathbb{R}\setminus\mathbb{Q}$, $n\in\mathbb{N}$. Obviously, $\Gamma_n\to\Gamma_0$ in $C^1(\mathbb{R})$. In particular, the densities $\Gamma_n^{\prime}(t)=\frac{f_n(t)}{A_n}$ and $\Gamma_0^{\prime}(t)=\frac{f_0(t)}{A_0}$ of invariant measures $\mu^{(n)}$ and $\mu^{(0)}$ converge uniformly. As for the interspike-interval distributions, we certainly have $\mu^{(n)}_{\textrm{ISI}}\implies\mu^{(0)}_{\textrm{ISI}}$. However, the assumption on the zero measure set $C_{\Psi_0}$ is not satisfied since $\Phi_0$ is a lift of the rotation and its  displacement $\Psi_0=\varrho_0$ is a constant function. Thus the set of critical points of $\Psi_0$ has full measure and indeed the distribution $\mu^{(0)}_{\textrm{ISI}}$ is degenerated to a point $\varrho_0$ (i.e. it is not absolutely continuous with respect to the Lebesque measure). Nevertheless, the distributions $\mu^{(n)}_{\textrm{ISI}}$ are absolutely continuous since the sets $C_{\Psi_n}$ of critical points of displacement functions $\Psi_n$ are countable (and the densities $\Delta^{(n)}(y)$ exist on the ground of Theorem \ref{wzor na gestosc}). Indeed: Note that $t\in\mathbb{R}$ is a critical point of $\Psi_n$ if and only if $\Phi_n^{\prime}(t)=1$. But $\Phi_n^{\prime}(t)=\frac{f_n(t)}{f_n(\Phi_n(t))}=1$ means that $\cos(2\pi n t)=\cos(2\pi n\Phi_n(t))$ which holds if and only if $\sin(\pi n(\Phi_n(t)+t))=0$ or $\sin(\pi n(\Phi_n(t)-t))=0$ as follows from the formula $\cos(\alpha)-\cos(\beta)=-2\sin(\frac{\alpha+\beta}{2})\sin(\frac{\alpha-\beta}{2})$. Further, $\sin(\pi n(\Phi_n(t)-t))=0$ only for $\Phi_n(t)-t=\frac{k}{n}$, where $k\in\mathbb{N}$, thus only for countably many choices of values $v_k=\frac{k}{n}$ of the displacement $\Phi_n(t)-t$. Now we will show that for every $n\in\mathbb{N}$ each value $v$ of the displacement function $\Psi_n(t)=\Phi_n(t)-t$ can be attained for at most countably many arguments $t\in\mathbb{R}$: Fix $v\in\mathbb{R}$. If $v=\Phi(t)-t$ for some $t\in \mathbb{R}$, we might assume that $v\in(0,1)$ (values of the displacement function are always contained is some interval of length not greater than 1, moreover $v=0$ or $v=1$ for some $t$ would imply that the firing phase map $\varphi_n$ has a fixed point at $x=\ue^{2\pi\imath t}$ which contradicts the irrationality of the rotation number $\varrho_n$). From definition of the firing map we have
\begin{equation}
\begin{split}
1 &=A_n(\Phi_n(t)-t)+\frac{B_n}{2\pi n}\sin(2\pi n\Phi_n(t))-\frac{B_n}{2\pi n}\sin(2\pi n t)\nonumber\\
\empty &=A_nv+\frac{B_n}{2\pi n}\sin(2\pi n(v+t))-\frac{B_n}{2\pi n}\sin(2\pi nt)
\end{split}
\end{equation}
This equality, after some calculations, leads to
\[
\frac{\pi n(1-A_nv)}{B_n}=\sin(\pi n v)\cos(\pi n(2t+v)).\] For $v\in(0,1)$ the above  equation is equivalent to
\[D_{n,v}=\cos(\widetilde{t}),\]
where $D_{n,v}=\frac{\pi n(1-A_nv)}{B_n(\sin(\pi n v))}$ and $\widetilde{t}=\pi n(2t+v)$. For fixed $n$ and $v$ $D_{n,v}$ is constant and this equality can be satisfied for at most countably many $\widetilde{t}\in\mathbb{R}$ and thus for every $n\in\mathbb{N}$ and $v\in\mathbb{R}$ also the entire equation $v=\Phi_n(t)-t$ is satisfied for at most countably many $t\in\mathbb{R}$. Now, since there are only countably many choices of values $v_k=\Phi_n(t)-t$ such that $\sin(\pi n(\Phi_n(t)-t))=0$ and each value $v_k$ is attained for at most countably many arguments $t\in\mathbb{R}$, $\sin(\pi n(\Phi_n(t)-t))=0$ for at most countably many $t\in\mathbb{R}$. Even easier one can justify that $\sin(\pi n(\Phi_n(t)+t))=0$ for at most countably many $t\in\mathbb{R}$, since $\Phi_n(t)+t$ is strictly increasing and thus it attains each value $v$ exactly once.  Thus for each $n\in \mathbb{N}$ the displacement functions $\Psi_n$ have countably many critical points. Theorem \ref{wzor na gestosc} yields now that the distributions $\mu^{(n)}_{\textrm{ISI}}$ have densities $\Delta^{(n)}$ with respect to the Lebesque measure. Nevertheless, the limiting (in terms of weak convergence of measures) distribution $\mu^{(0)}_{\textrm{ISI}}$ does not have the density. \\

\noindent{\textbf{Example 4.}} Let us consider Perfect-Integrator Model: $\dot{x}=f_n(t)$, where \[f_n(t)=A+B\cos(2\pi nt), \ n\in\mathbb{N}.\] If $A>B>0$ then $f_n(t)>0$ for every $n\in\mathbb{N}$ and $t\in\mathbb{R}$. In this case the firing maps $\Phi_n:\mathbb{R}\to\mathbb{R}$ induced by the equations $\dot{x}=f_n(t)$ are lifts of orientation preserving circle homeomorphisms $\varphi_n:S^1\to S^1$. Moreover, each $\Phi_n$ is conjugated with the lift $R_{\varrho_n}(t)=t+\frac{1}{A}$ of the rotation by $\varrho_n=1/\int_{0}^{1}A+B\cos(2\pi n t)\mathrm{dt}=\frac{1}{A}$  via
\[
\Gamma_n(t)=\frac{\int_{0}^{t}A+B\cos(2\pi nu)\mathrm{du}}{A}=t+\frac{B}{2\pi n A}\sin(2\pi n t)
\]
In particular, all the rotation numbers $\varrho_n$ of $\Phi_n$ are the same and can be set rational or irrational with arbitrary diophantine properties (it depends only on the choice of $A$).

Since $\Gamma_n\to \textrm{Id}$ uniformly (i.e. in $C^0(\mathbb{R})$), also $\Gamma_n^{-1}\to \textrm{Id}$ uniformly. Consequently,
\[
\Phi_n(t)=\Gamma_n^{-1}(\Gamma_n(t)+\varrho_n)=\Gamma_n^{-1}(t+\frac{B\sin(2\pi n t)}{2\pi n A}+\frac{1}{A})\to t+\frac{1}{A}
\]
and thus also $\Phi_n\to \Phi_0$ uniformly with $\Phi_0(t)=t+\frac{1}{A}$ being simply a lift of the rotation by $\frac{1}{A}$. Note that $\Phi_0(t)$ can be seen as a firing map induced by the equation $\dot{x}=f_0(t)$ with $f_0(t)=A$. However, it is not true that $f_n\to f_0$ uniformly or even pointwise since for example for $t=\frac{1}{2}$ the sequence $f_n(t)=A+B\cos(2\pi n t)$ does not converge at all.

From the uniform convergence $\Phi_n\to\Phi_0$ we have the weak convergence $\mu_n\implies \mu_0$ of the corresponding unique invariant probability measures:
\[
\mu_{n}(V)=\frac{\int_V f_n(u)\mathrm{du}}{A}, \quad V\subset [0,1]
\]
where $\mu_0=\Lambda$ is the Lebesque measure on $[0,1]$ (being the invariant measure of $\Phi_0$). From the formula for $\mu_n$ we see that the invariant measures $\mu_{n}$ have densities
\[
\widetilde{f}_n=\frac{f_n}{A}
\] and the measure $\mu_0$ has a density
\[
\widetilde{f}_0=\frac{f_0}{A}\equiv 1
\]
However, $\widetilde{f_n}\nrightarrow \widetilde{f}_0$, similarly as $f_n \nrightarrow f_0$. In particular, this shows that the weak convergence of measures does not imply (even pointwise) convergence  of the corresponding continuous density functions. Thus the sequence of conjugacies $\Gamma_n$ does not converge in $C^1(\mathbb{R})$ but only in $C^0(\mathbb{R})$. \\
\subsection{Empirical approximation of the interspike-interval distribution $\mu_{\textrm{ISI}}$}
 Virtually we are able to calculate  only the empirical interspike-interval distribution, i.e. the distribution derived by counting interspike-intervals along a particular trajectory (a run of a system). In case of the rational firing rate necessarily there are periodic orbits (and usually also non-periodic but these are attracted to the periodic ones) and although all the periodic orbits have the same period, the (finite) sequences $\textrm{ISI}_n(t_0)$ derived along the orbit of each periodic point $t_0$ might consist of different elements unless the system induces a rigid rotation. However, in case of the irrational rotation number we have the unique invariant ergodic measure that gives the distribution of orbits phases. Thus $\mu_{\textrm{ISI}}$ is also well-defined and the empirical distribution of interspike-intervals derived for any trajectory will  well approximate $\mu_{\textrm{ISI}}$, provided that the trajectory is long enough. However, basically when we do numerical computations, we do not work with irrational rotation numbers, but the rational ones which are close to them. We will see in what meaning the empirically derived interspike-interval distribution for an arbitrarily chosen initial condition $(t,0)$ of a system with rational firing rate, being ``close'' enough to our entire system with irrational firing rate, approximates the desired distribution $\mu_{\textrm{ISI}}$ of the ``irrational'' (ergodic) system.

We have to  define the empirical interspike-interval distribution formally:
\begin{definition}
Let $\Phi$ be the firing map arising from the IF system $\dot{x}=F(t,x)$. Choose the initial condition $(t,0)$ ($x_r=0$). Then the empirical interspike-interval distribution for the run of length $n$ (i.e. having $n$-spikes) starting at $(t,0)$ equals
\begin{displaymath}
\omega_{n,t}=\frac{1}{n}\sum_{i=0}^{n-1}\delta_{\textrm{ISI}_{i}(t)},
\end{displaymath}
where $\delta_{\textrm{ISI}_{i}(t)}$ is a Dirac delta centered at the point $\textrm{ISI}_i(t)=\Phi^{i+1}(t)-\Phi^{i}(t)$.
\end{definition}
Thus $\omega_{n,t}(A)=\frac{1}{n}\sharp\{0\leq i\leq n-1: \ \Phi^{i+1}(t)-\Phi^{i}(t)\in A\}$, $A\subset \mathbb{R}$.

 Note that if $\widetilde{\varphi}$ with rotation number $\widetilde{\varrho}$ is close in $C^0(S^1)$-metric to $\varphi$ with irrational rotation number $\varrho$, then the rotation numbers $\widetilde{\varrho}$ and $\varrho$ are also close due to the continuity of the rotation number in $C^0(S^1)$.

In order to measure the distance between interspike-interval distribution we introduce the notion of the Fortet-Mourier metric:
\begin{definition} Let $\mu$ and $\nu$ be the two Borel probability measures on a measurable space $(\Omega, \mathcal{F})$, where $\Omega$ is a compact metric space. Then the distance between the measures $\mu$ and $\nu$ is defined as
\begin{displaymath}
d_F(\mu,\nu):=\sup\{ \vert\int\limits{\Omega} f\; d\mu - \int\limits{\Omega} f\; d\nu\vert : \ f \ \textrm{is $1$-Lipschitz}\}.
\end{displaymath}
\end{definition}
Using Theorem 3.15 in \cite{wmjs4} we formulate the following:
\begin{proposition}\label{przyblizaniedlalif}
Consider the integrate-and-fire systems $\dot{x}=-\sigma x +f_{\theta_1}(t)$ and $\dot{x}=-\sigma x +f_{\theta_2}(t)$, where $f_{\theta_i}\in L^{\infty}_{\textrm{loc}}(\mathbb{R})$, periodic with period $1$ and $f_{\theta_1}(t)-\sigma>\varsigma>0$. By $\Phi_{\theta_1}$ and $\Phi_{\theta_2}$ denote the firing maps emerging from the corresponding systems. Suppose that the rotation number associated with $\Phi_{\theta_1}$ is irrational.

For any $\varepsilon>0$ there exists a neighbourhood $\mathcal{U}$ of $f_{\theta_1}$ in $L^{\infty}_{\textrm{loc}}(\mathbb{R})$-topology such that if $f_{\theta_2}\in\mathcal{U}$, then for every initial condition $(t,0)$ we have:
\begin{equation}
d_F(\lim_{n\to\infty}\omega^{(\theta_2)}_{n,t}, \mu^{(\theta_1)}_{\textrm{ISI}})<\varepsilon,
\end{equation}
where $\omega^{(\theta_2)}_{n,t}$ is the empirical interspike-interval distribution for the run of the system $\dot{x}=-\sigma x + f_{\theta_2}(t)$ starting from $(t,0)$ and $\mu^{(\theta_1)}_{\textrm{ISI}}$ is the interspike-interval distribution for $\dot{x}=-\sigma x + f_{\theta_1}(t)$ with respect to its invariant measure $\mu^{(\theta^1)}$.
\end{proposition}
\noindent{\bf Proof.} The proof relies again of the fact that the mapping $f\mapsto\Phi$ is continuous from $\textrm{ess sup}$-topology into $C^0(\mathbb{R})$. Then the statement follows immediately  from  Theorem 3.15 in \cite{wmjs4}.\hbx\\

\begin{figure}[h!]
\begin{subfigure}[h]{0.48\textwidth}
                \centering
                \includegraphics[width=0.49\textwidth]{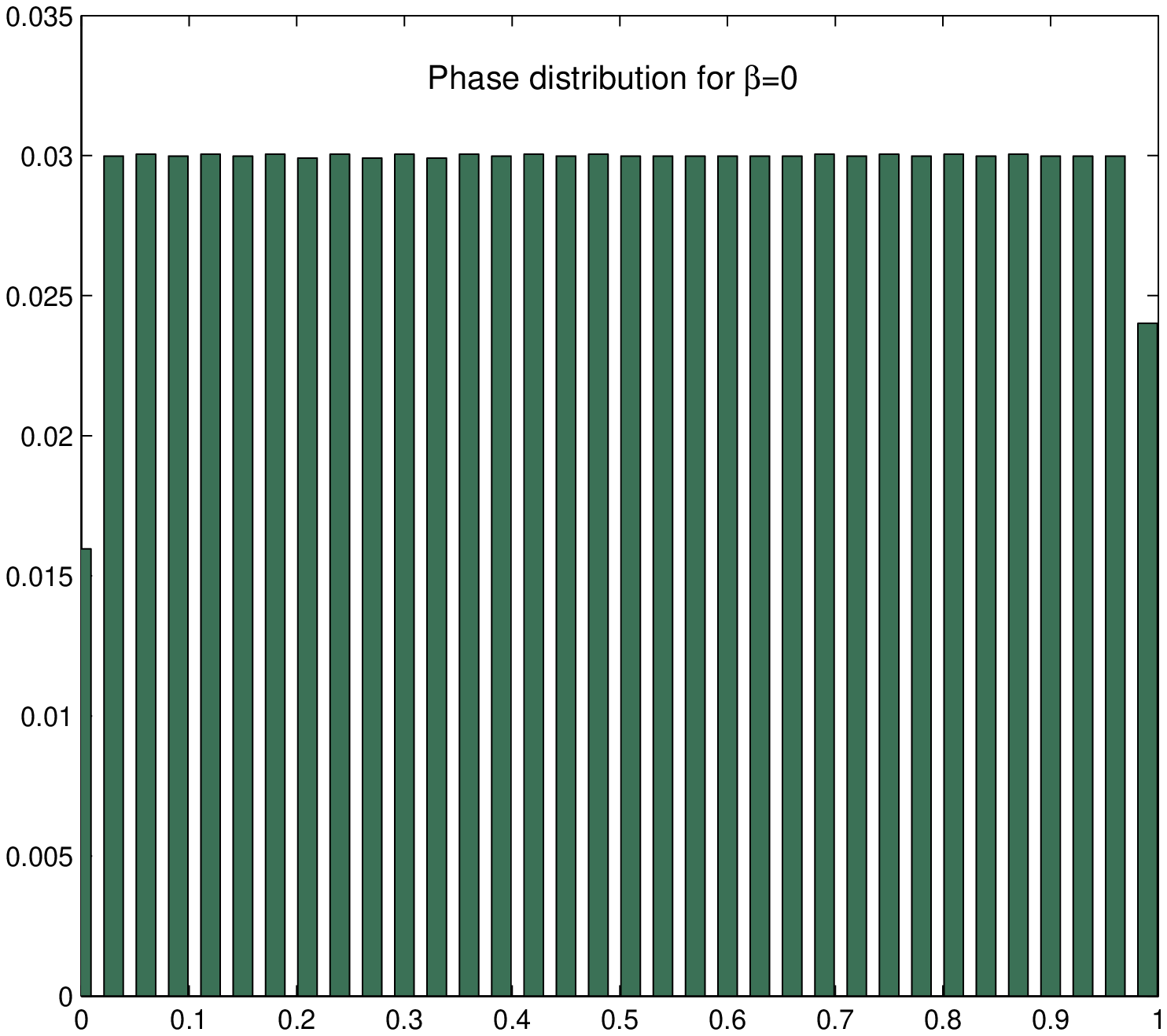}
                \includegraphics[width=0.49\textwidth]{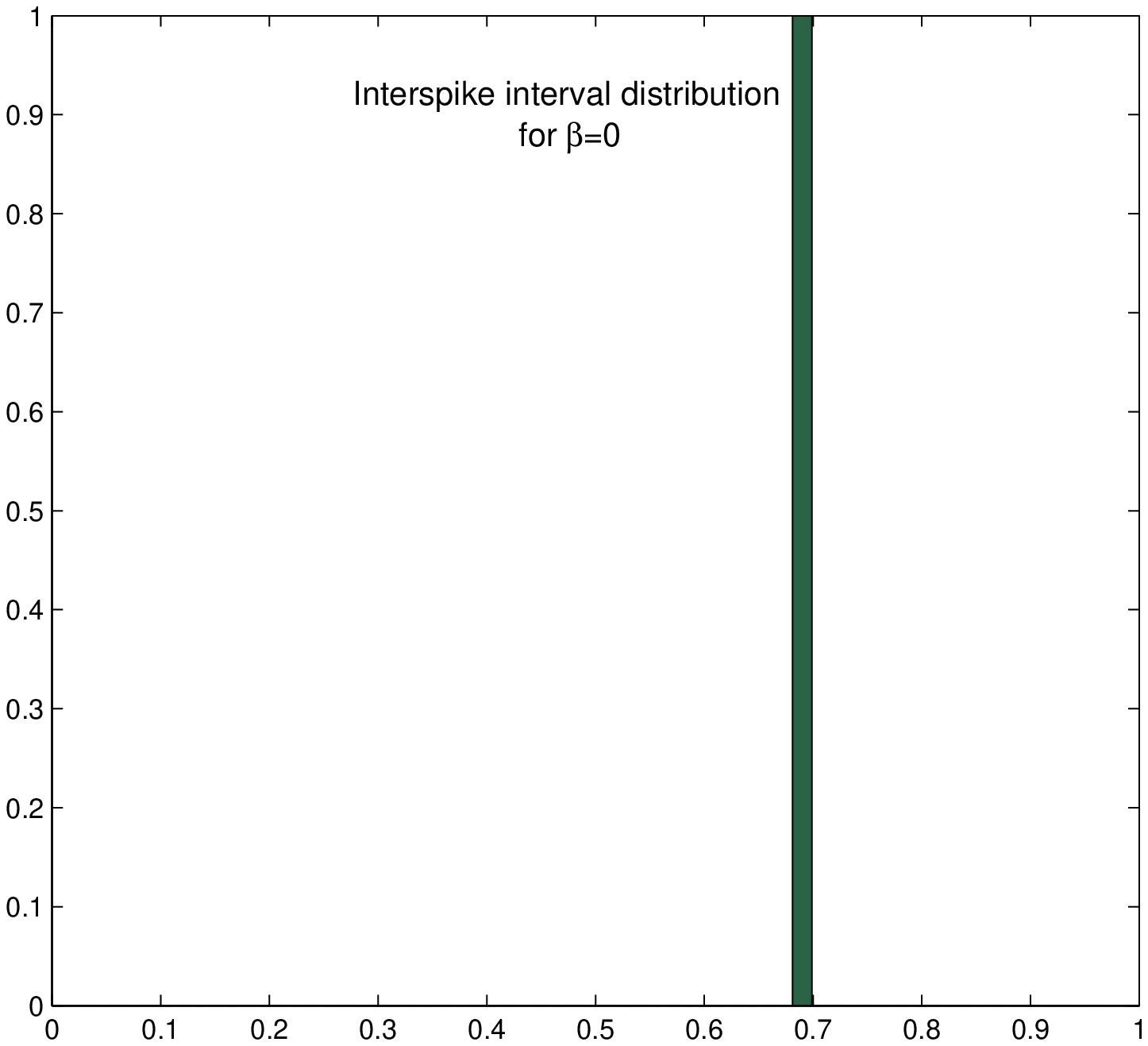}
                \caption{Firing phases (left) and interspike-interval (right) distribution for $\beta=0$.}
                \label{fig:beta0}
        \end{subfigure}
\hspace{1.5cm}
        \begin{subfigure}[h]{0.48\textwidth}
                \centering
                \includegraphics[width=0.49\textwidth]{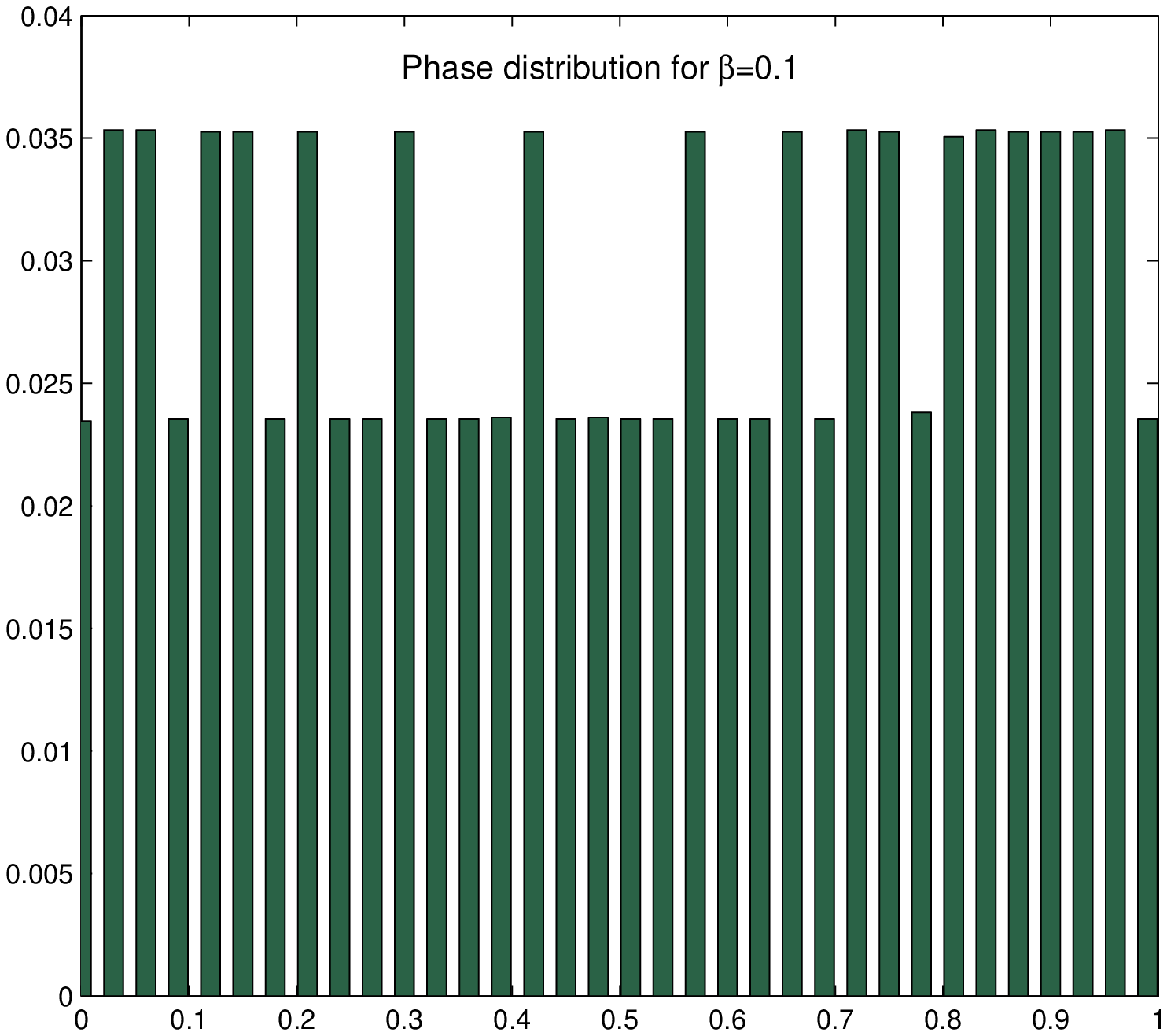}
                \includegraphics[width=0.49\textwidth]{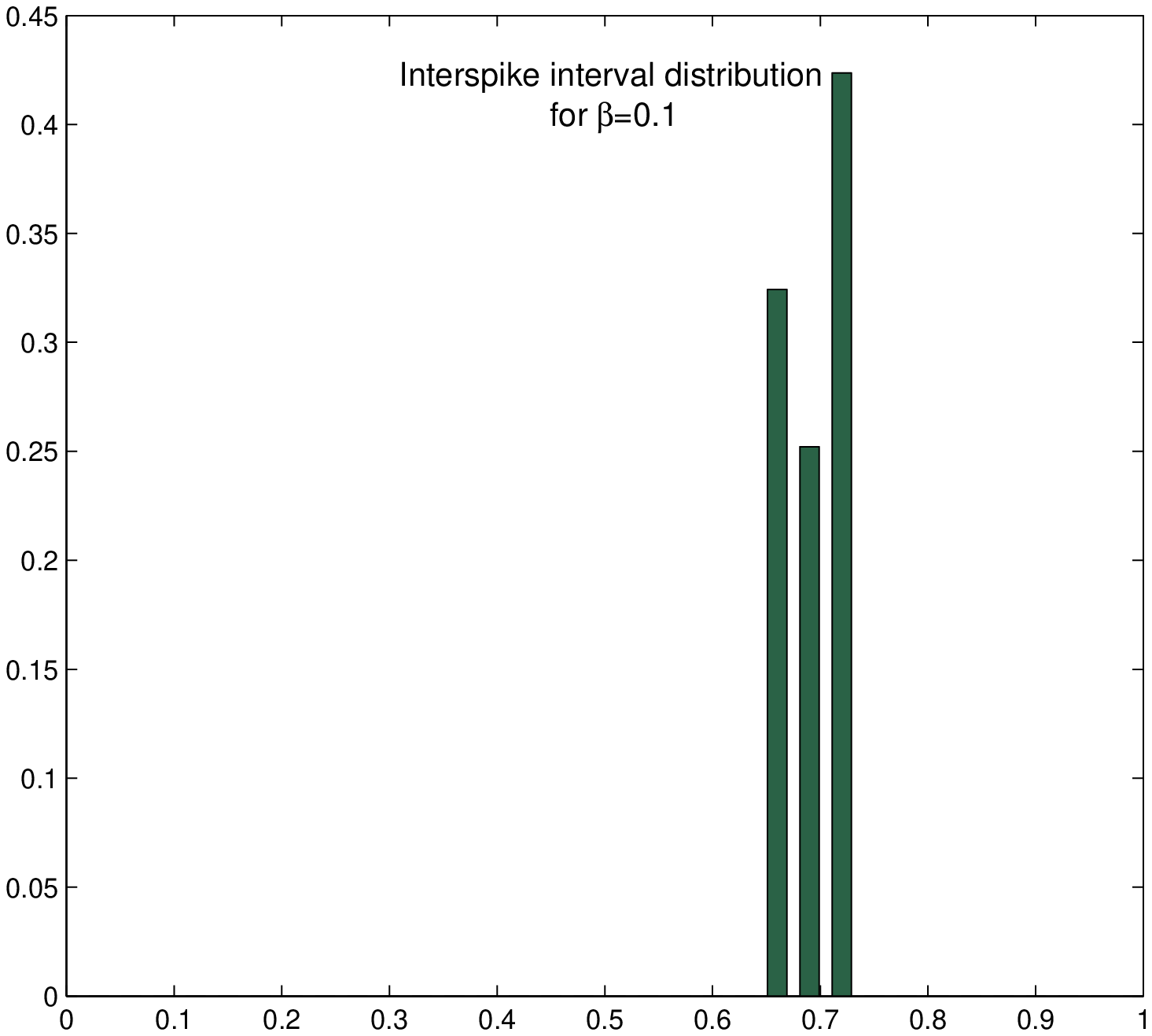}
                \caption{Firing phases (left) and interspike-interval (right) distribution for $\beta\nolinebreak[10]=\nolinebreak[10]0.1$.}
                \label{fig:beta01}
        \end{subfigure}

\vspace{1cm}
        \begin{subfigure}[h]{0.48\textwidth}
                \centering
                \includegraphics[width=0.49\textwidth]{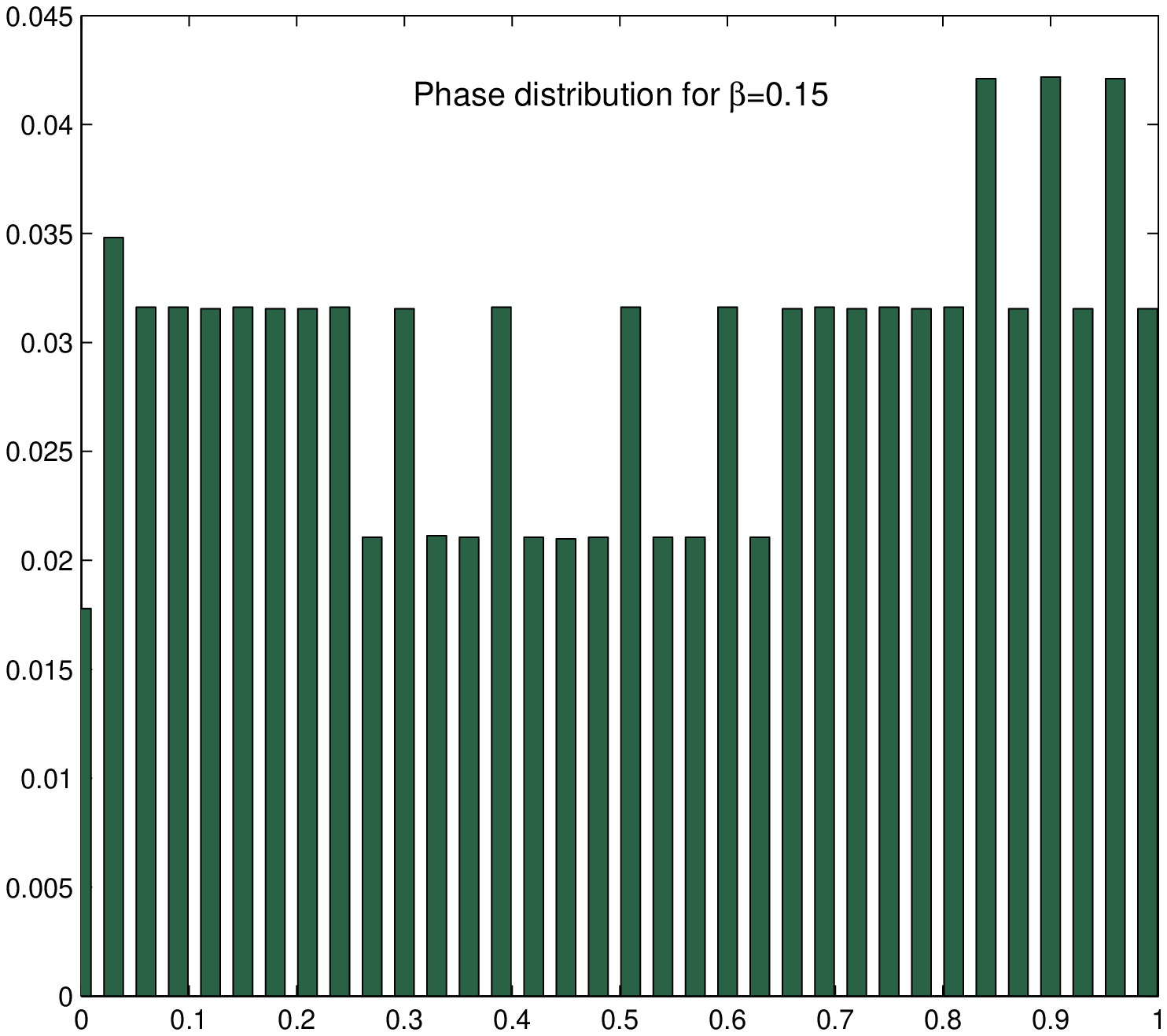}
                \includegraphics[width=0.49\textwidth]{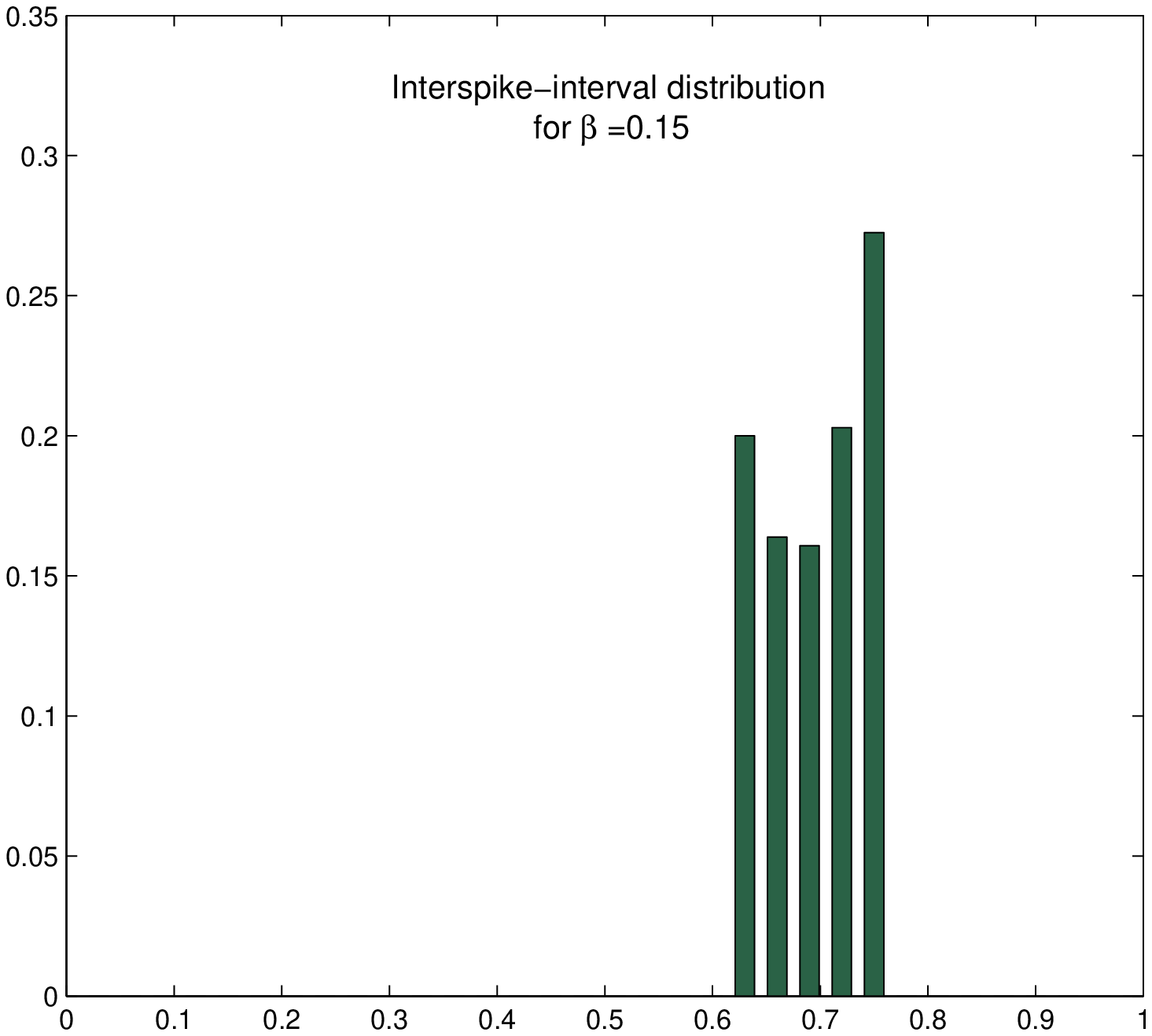}
                \caption{Firing phases (left) and interspike-interval (right) distribution for $\beta\nolinebreak[10]=\nolinebreak[10]0.15$.}
                \label{fig:beta015}
        \end{subfigure}
\hspace{1.5cm}
        \begin{subfigure}[h]{0.48\textwidth}
                \centering
                \includegraphics[width=0.49\textwidth]{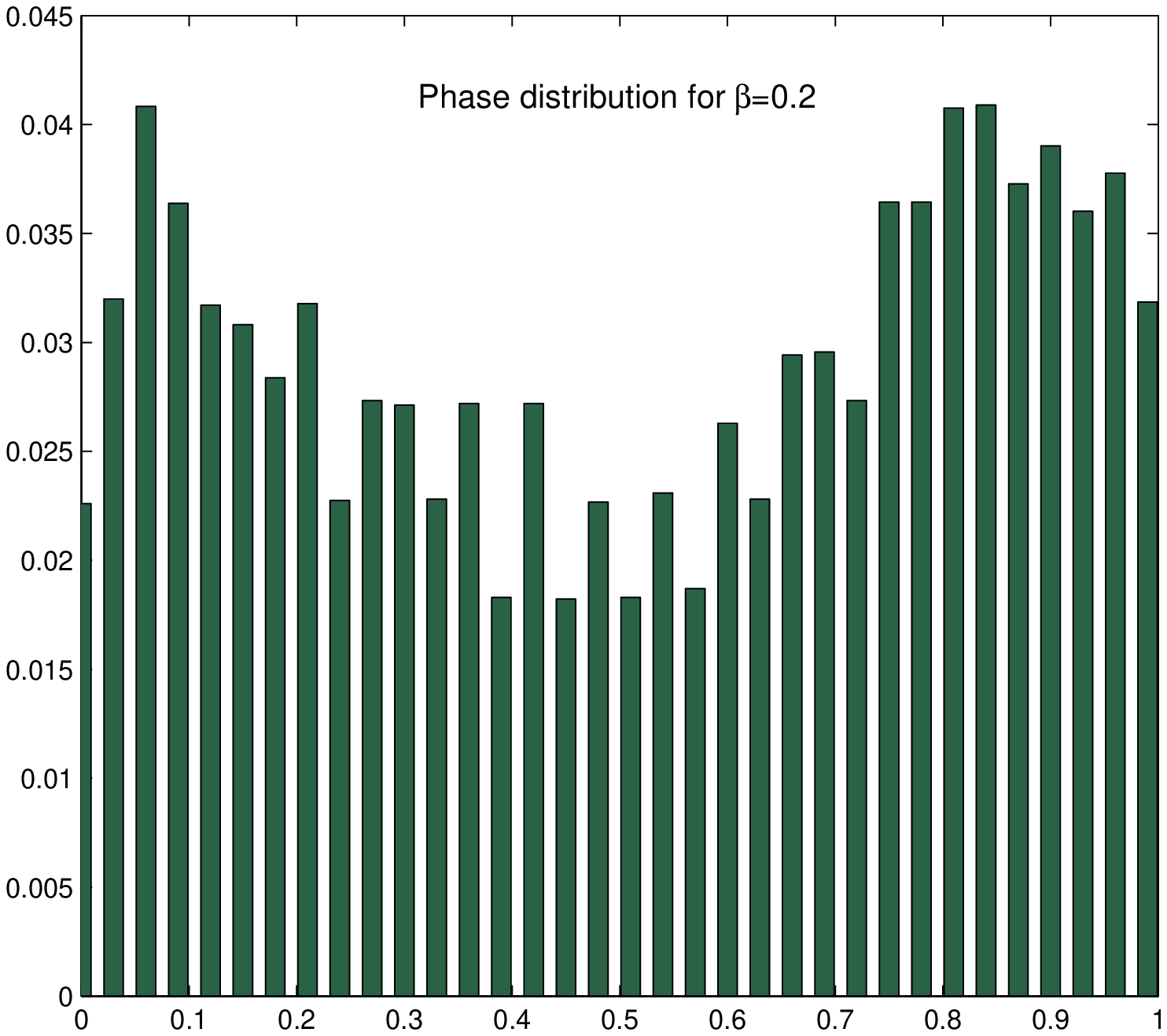}
                \includegraphics[width=0.49\textwidth]{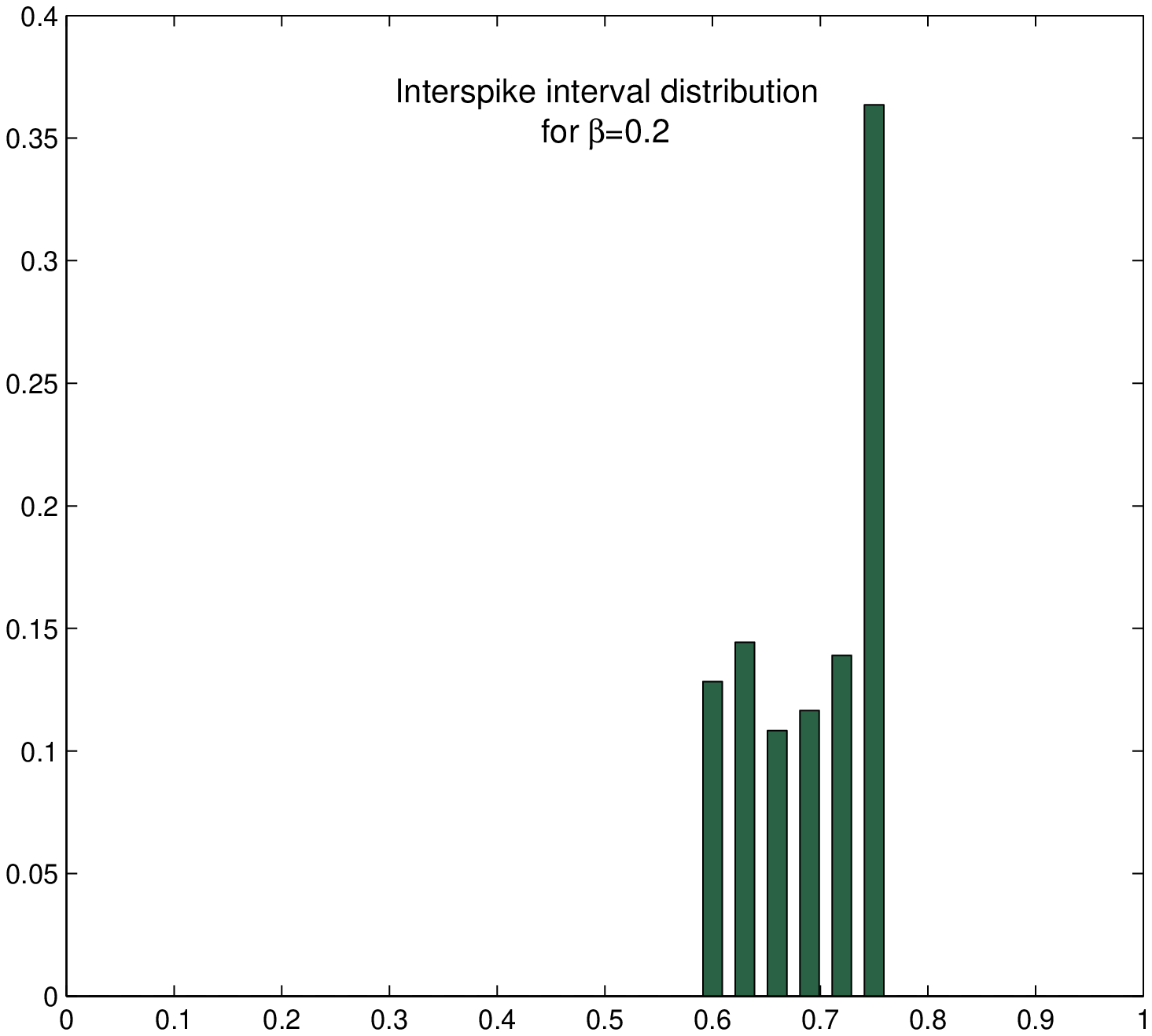}
                \caption{Firing phases (left) and interspike-interval (right) distribution for $\beta\nolinebreak[10]=\nolinebreak[10]0.2$.}
                \label{fig:beta02}
        \end{subfigure}

\vspace{1cm}
        \begin{subfigure}[h]{0.48\textwidth}
                \centering
                \includegraphics[width=0.49\textwidth]{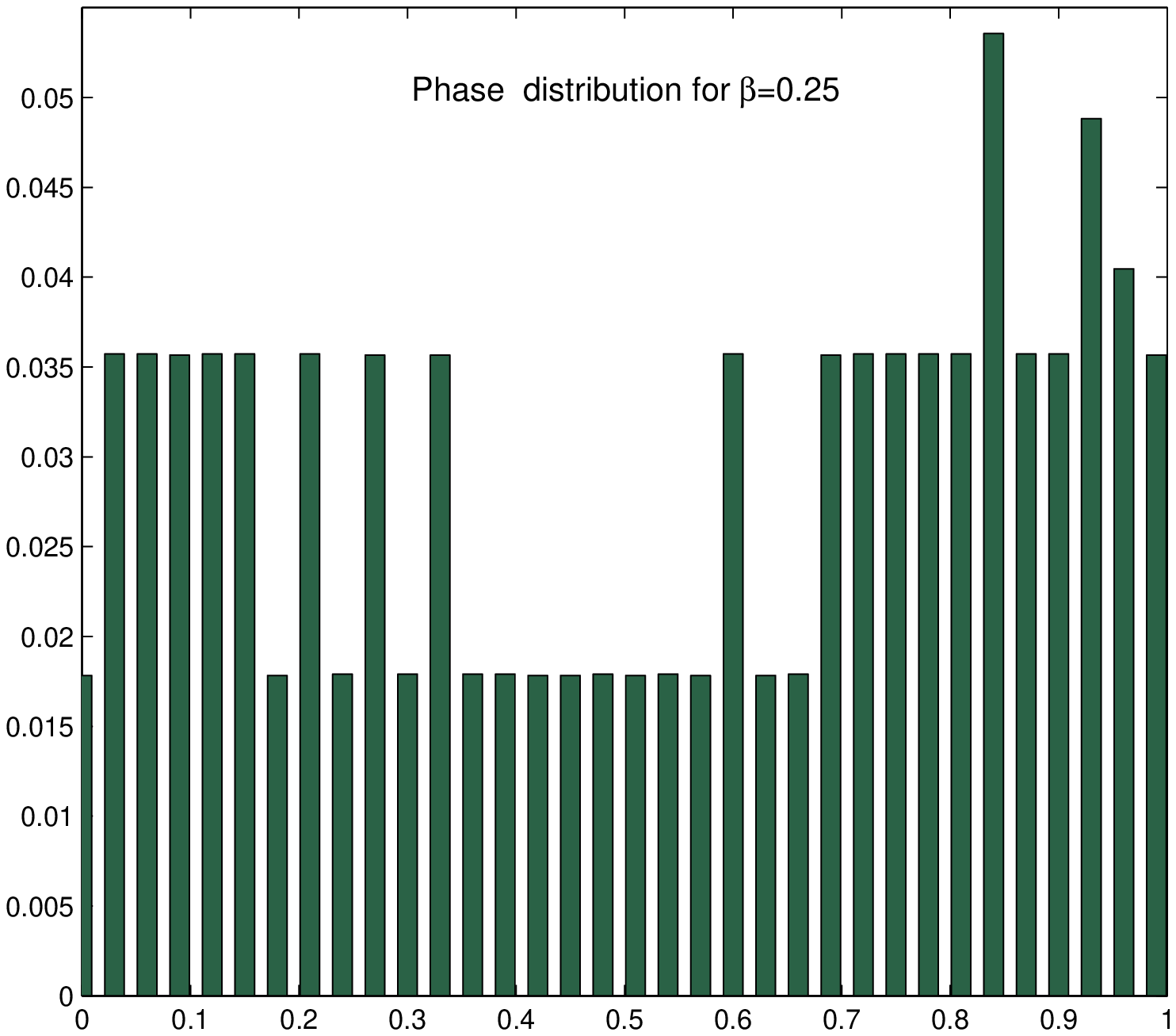}
                \includegraphics[width=0.49\textwidth]{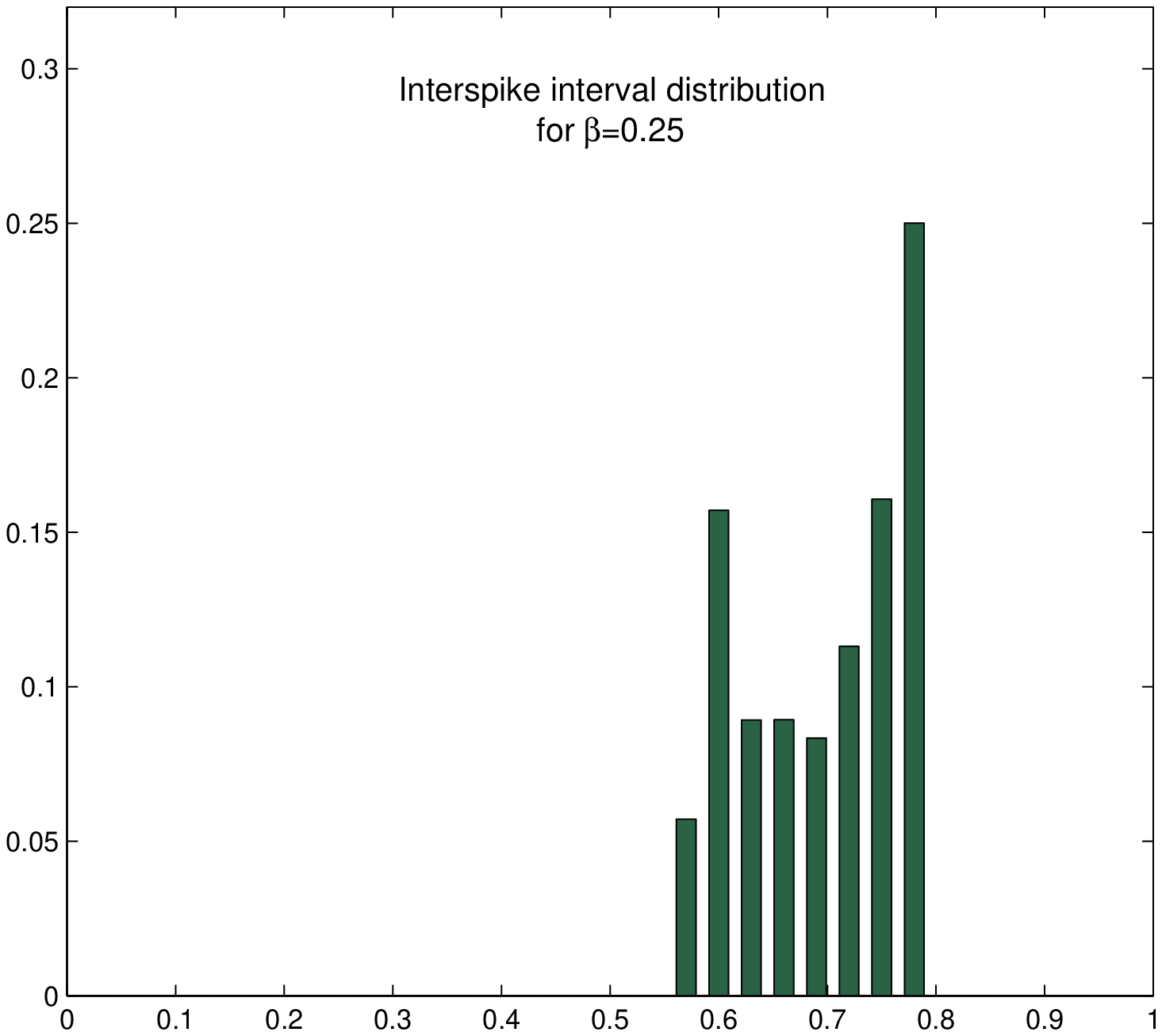}
                \caption{Firing phases (left) and interspike-interval (right) distribution for $\beta\nolinebreak[10]=\nolinebreak[10]0.25$.}
                \label{fig:beta025}
        \end{subfigure}

        \caption[Distribution of firing phases and interspike-intervals for the model $\dot{x}=-\sigma x + 2(1+\beta\cos(2\pi t))$ and $\beta=0,0.1, 0.15,0.2, 0.25$.]{Distribution of firing phases  and interspike-intervals  for the model $\dot{x}=-\sigma x + 2(1+\beta\cos(2\pi t))$ and $\beta=0,0.1, 0.15,0.2, 0.25$.}\label{fig:beta1}
\end{figure}

As the convergence under Fortet-Mourier metric implies weak convergence of measures (\cite{probmetric}) we conclude:
\begin{corollary}
Under the notation as in Proposition \ref{przyblizaniedlalif}, for every $t\in\mathbb{R}$ we have
\begin{displaymath}
\lim_{n\to\infty}\widetilde{\omega}_{n,t}^{(\theta_2)}\implies \mu_{\textrm{ISI}}^{(\theta_1)}.
\end{displaymath}
\end{corollary}
The above result can be illustrated by the numerical example:\\
\noindent{\textbf{Example 5.}} We investigate the system $\dot{x}=-x +2(1+\beta\cos(2\pi t))$ (the computations were done in Matlab). Notice that the analogous example was considered in a classical paper \cite{keener1}, but the authors gave no explanation of the behaviour of interspike-intervals histograms under a  small change of a parameter. Our results allow us to make theoretical predicates of what actually we can expect for the interspike-interval distribution when the parameter $\beta$ varies. We easily obtain that for $0\leq\beta<0.5$ the firing map $\Phi:\mathbb{R}\to\mathbb{R}$ is a lift of a circle homeomorphism.  The results of numerical simulations for $\beta=0,0.1,0.15,0.2$ and $0.25$ are presented in Figure \ref{fig:beta1}. All the simulations were started from the initial condition $(0,0)$.

When $\beta=0$ the system is forced by the constant input which induces the rotation by an irrational angle. Indeed, by solving the equation and direct computation we obtain that  $\textrm{ISI}_n(t)=\varrho=\ln(0.5)\approx 0.6931$ for every $n$ and $t\in\mathbb{R}$. This is reflected in Figure \ref{fig:beta0}: the firing phases are distributed uniformly in $[0,1]$ and the iterspike-interval distribution is simply a Dirac delta at $\textrm{ISI}^{(0)}\approx0.6931$. Thus for $\beta=0$ we are dealing with the irrational rotation. When we slightly change the parameter $\beta$ (Figure \ref{fig:beta01}-\ref{fig:beta025}, we observe that both the distribution of phases and of interspike-intervals change continuously as we anticipate from the fact that the corresponding distributions are close in $d_F$ metric, since the firing maps are close in $C^0(\mathbb{R})$ metric. In particular the distribution of interspike-intervals is concentrated in the interval around the value of $\textrm{ISI}^{(0)}$ and the distribution practically fills this whole interval, which is consistent with Proposition \ref{gestosc w przedziale dla lif} as in our case the input function is smooth.

\begin{figure}[h!]
        \begin{subfigure}[h]{0.48\textwidth}
                \centering
\includegraphics[width=0.49\textwidth]{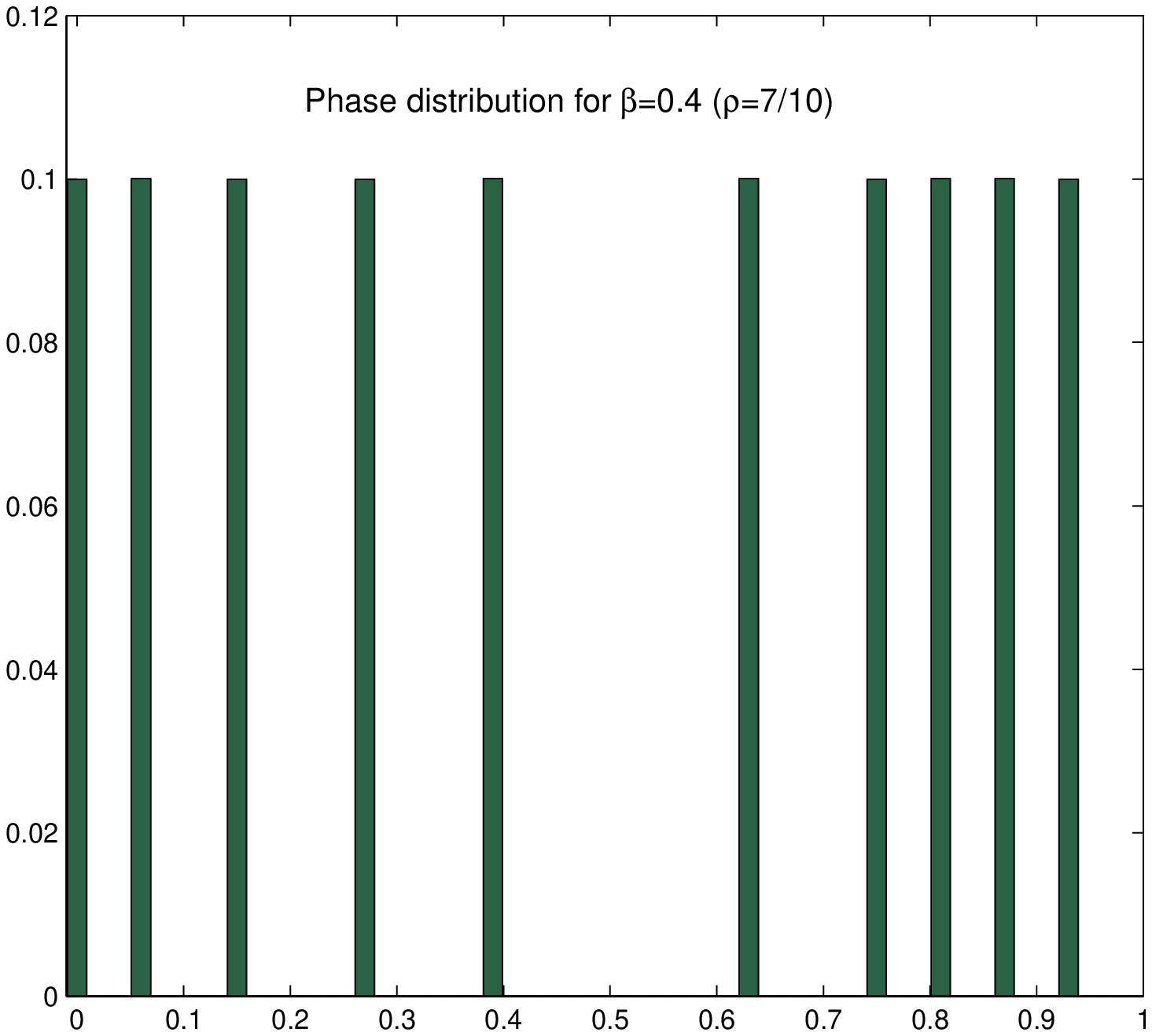}
\includegraphics[width=0.49\textwidth]{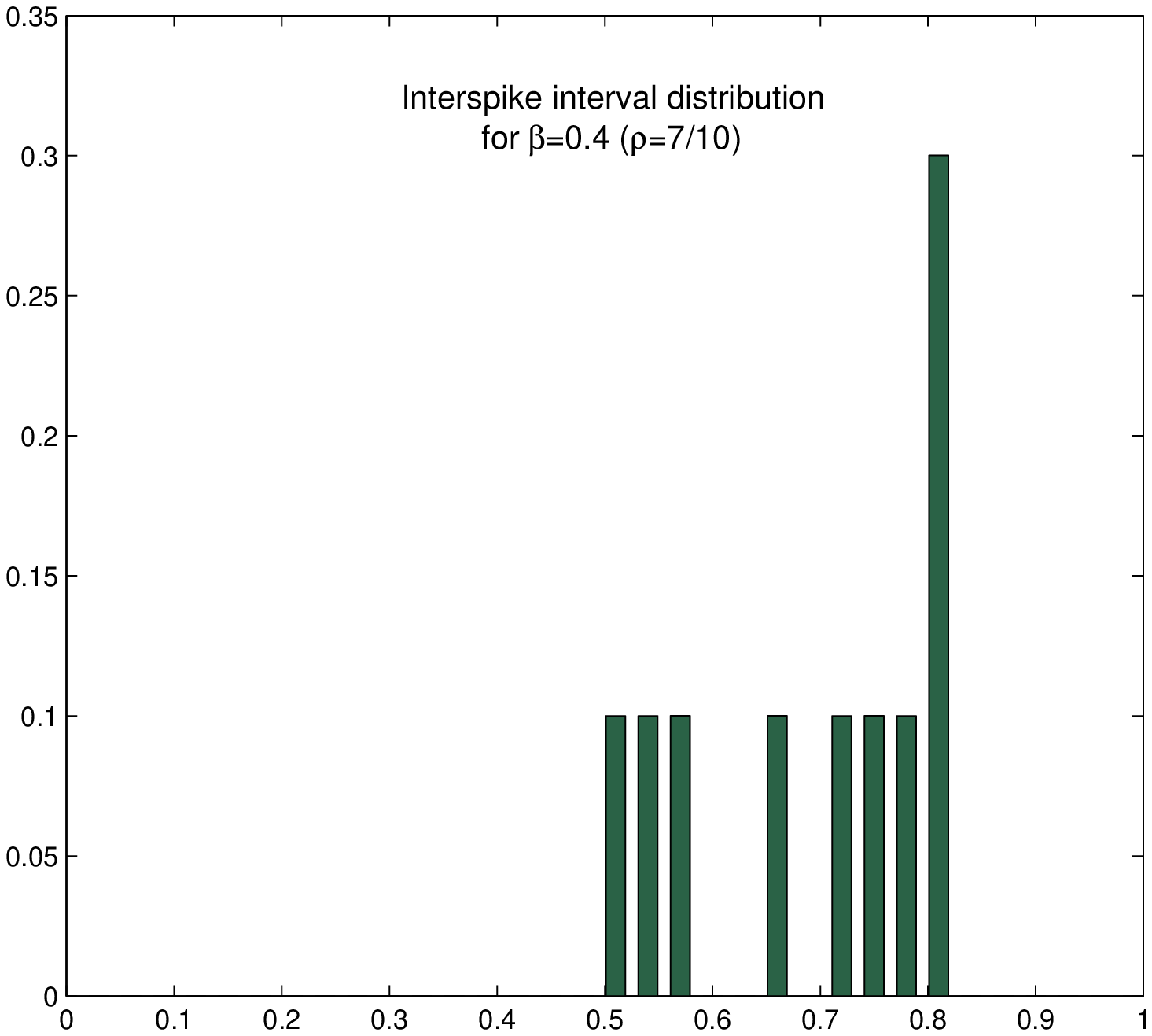}
 \caption{Firing phases (left) and interspike-interval (right) distribution for $\beta\nolinebreak[10]=\nolinebreak[10]0.4$.}
\label{fig:beta2_04}
 \end{subfigure}
 \hspace{1.5cm}
        \begin{subfigure}[h]{0.48\textwidth}
                \centering
                \includegraphics[width=0.49\textwidth]{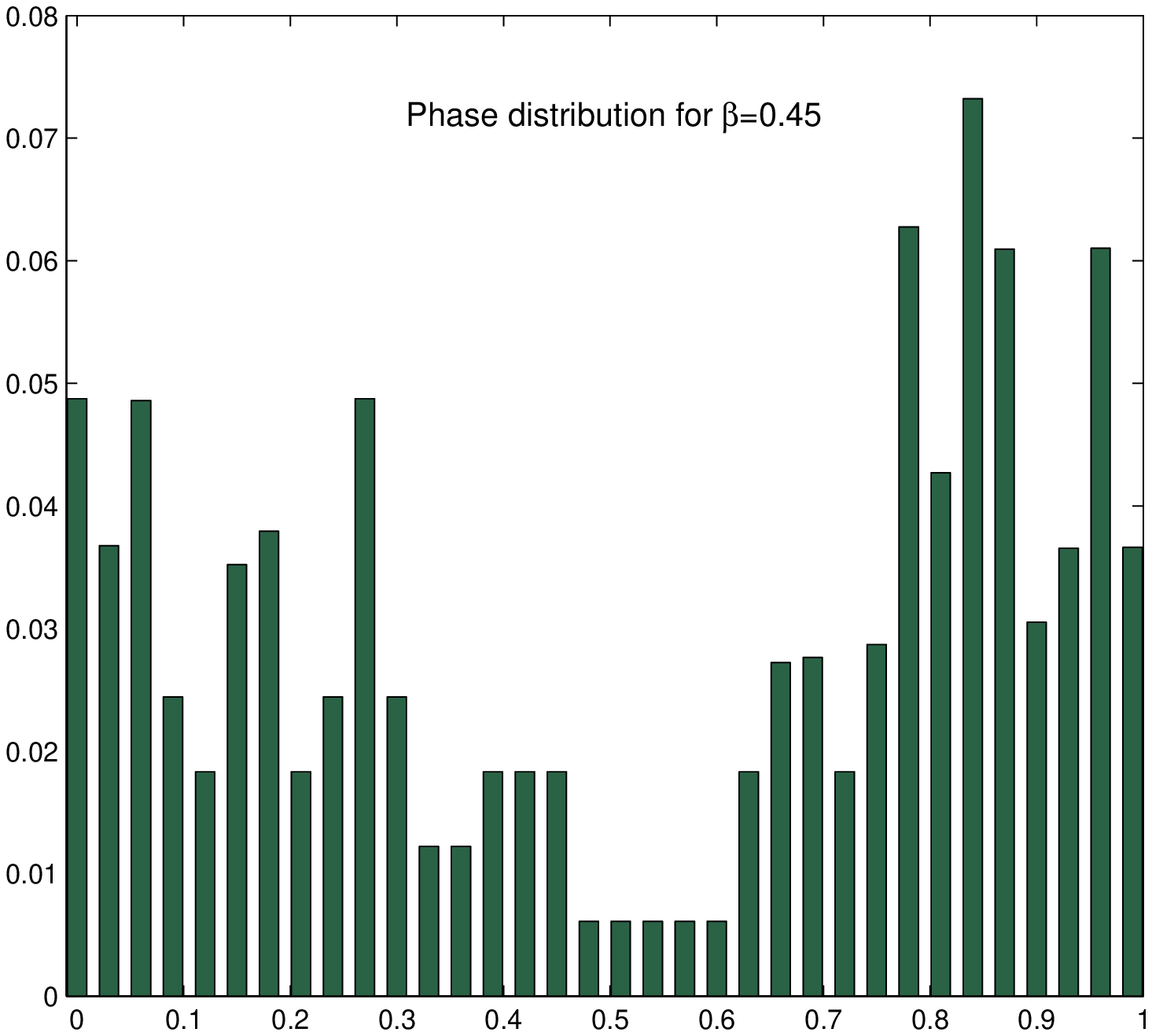}
                \includegraphics[width=0.49\textwidth]{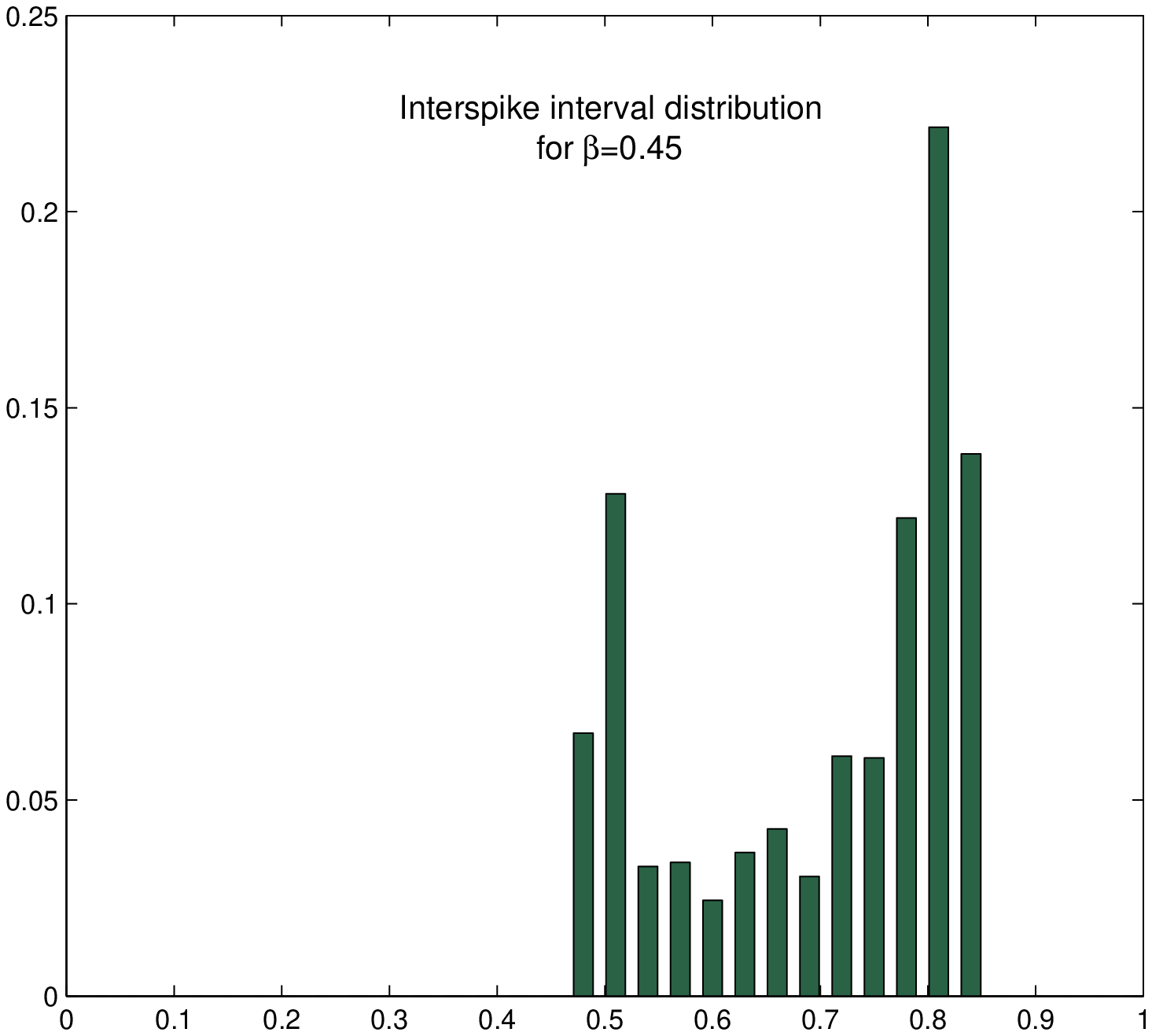}
                \caption{Firing phases (left) and interspike-interval (right) distribution for $\beta\nolinebreak[10]=\nolinebreak[10]0.45$.}
                \label{fig:beta2_045}
        \end{subfigure}

\vspace{1cm}
        \begin{subfigure}[h]{0.48\textwidth}
                \centering
                \includegraphics[width=0.49\textwidth]{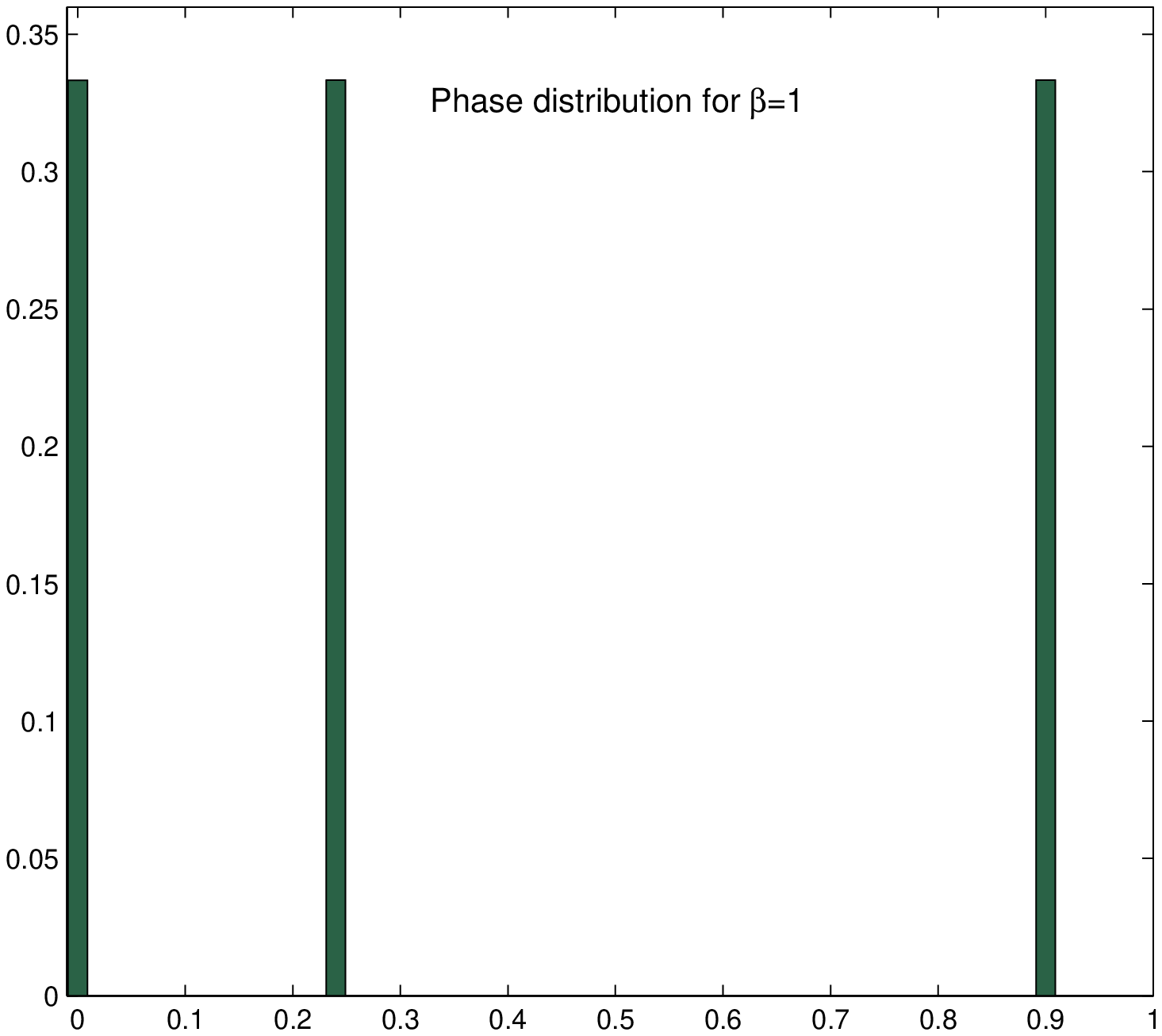}
                \includegraphics[width=0.49\textwidth]{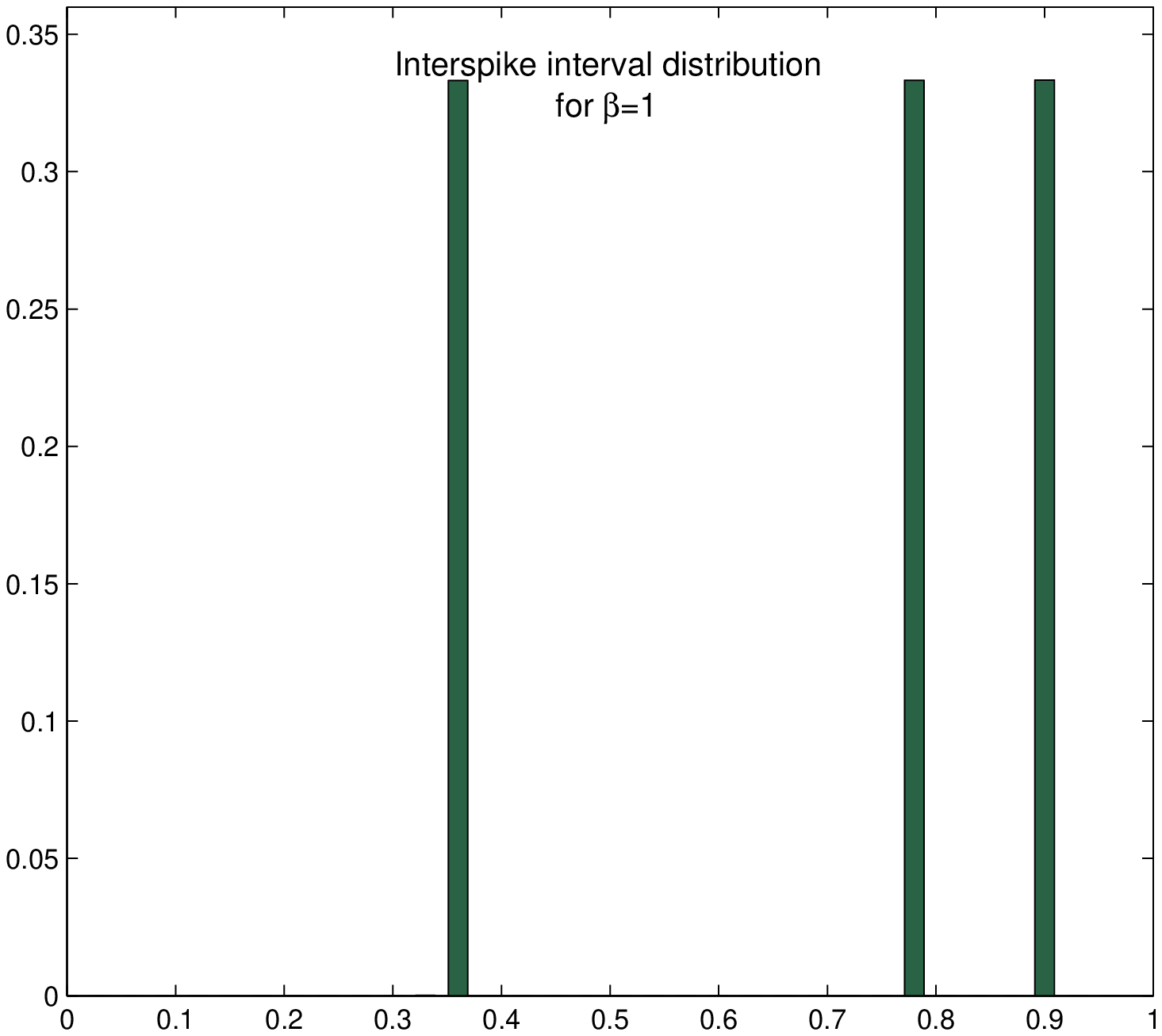}
                \caption{Firing phases (left) and interspike-interval (right) distribution for $\beta\nolinebreak[10]=\nolinebreak[10]1$.}
                \label{fig:beta2_1}
        \end{subfigure}
        \hspace{1.5cm}
        \begin{subfigure}[h]{0.48\textwidth}
                \centering
                \includegraphics[width=0.49\textwidth]{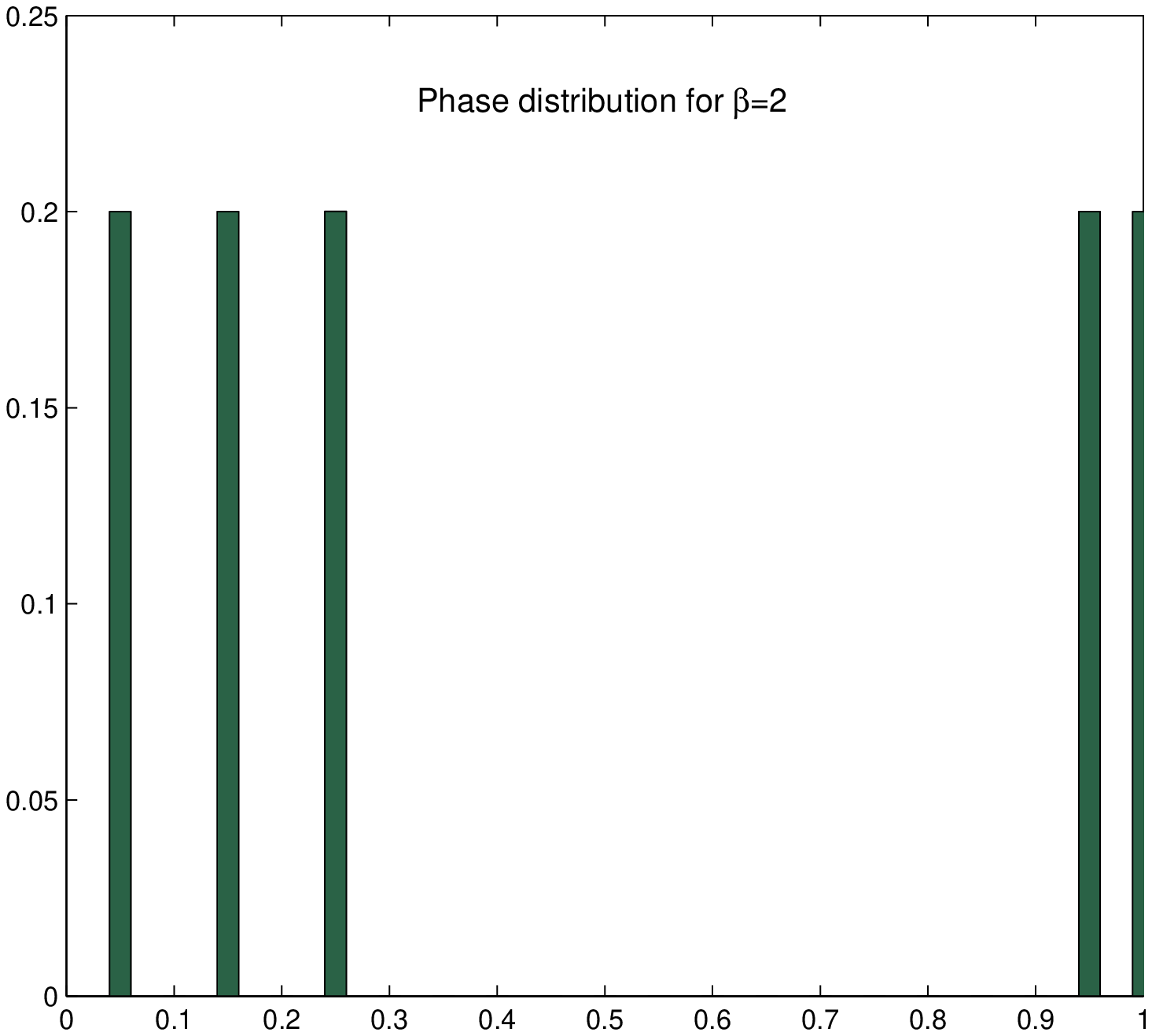}
                \includegraphics[width=0.49\textwidth]{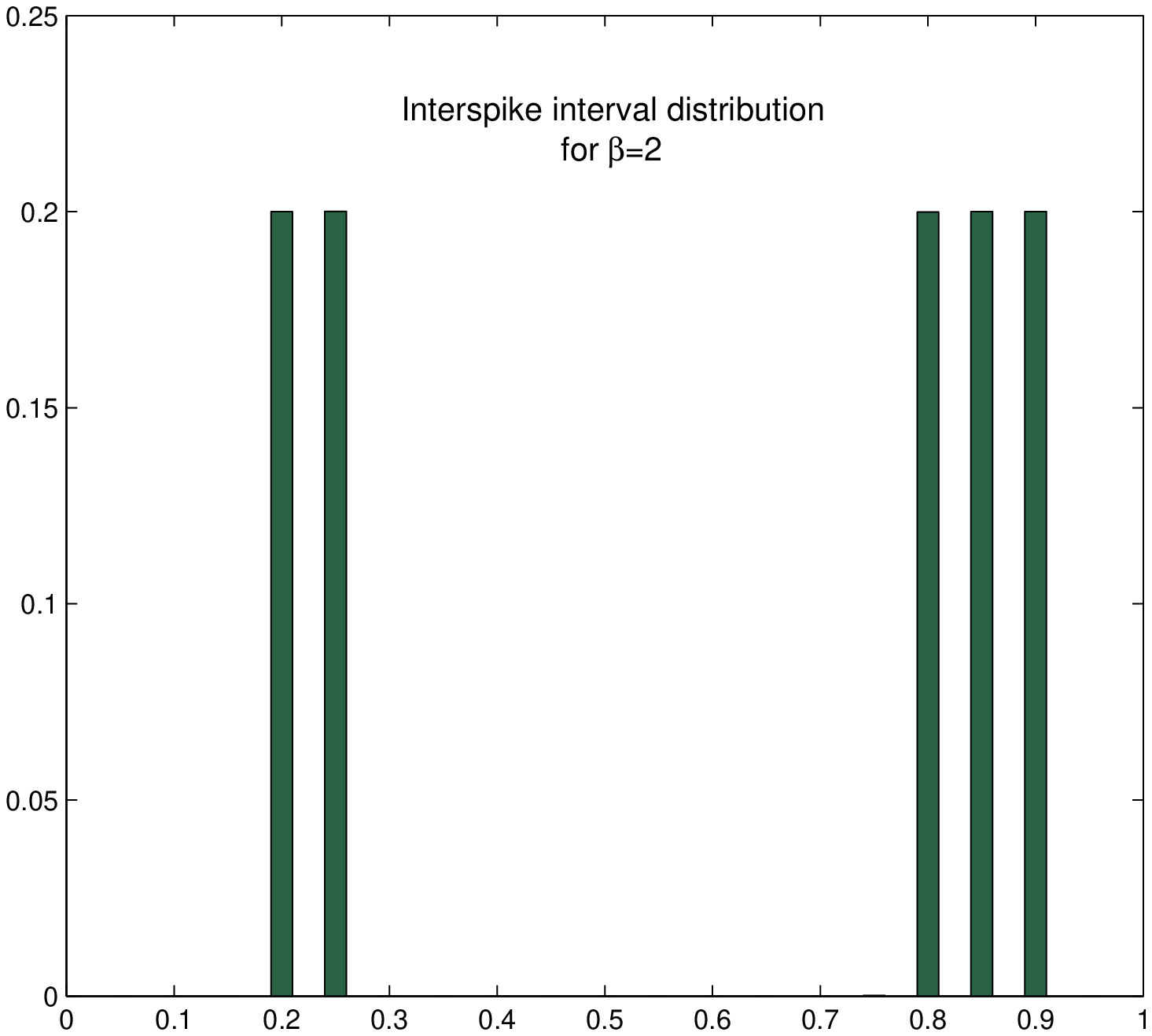}
                \caption{Firing phases (left) and interspike-interval (right) distribution for $\beta\nolinebreak[10]=\nolinebreak[10]2$.}
                \label{fig:beta2_2}
        \end{subfigure}

        \caption[Distribution of firing phases and interspike-intervals for the model $\dot{x}=-\sigma x + 2(1+\beta\cos(2\pi t))$ and $\beta=0.4, 0.45, 1, 2$.]{Distribution of firing phases  and interspike-intervals  for the model $\dot{x}=-\sigma x + 2(1+\beta\cos(2\pi t))$ and $\beta=0.4, 0.45, 1, 2$.}\label{fig:beta2}
\end{figure}

We also have checked what happens for greater values of $\beta$ and the results are presented in Figure \ref{fig:beta2}. When $\beta=0.4$, both firing phases and interspike-intervals admit ten distinct values, which suggest that there is a periodic orbit of period $10$ and indeed, the rotation number was computed as $\varrho=7/10$. When the parameter changes to $\beta=0.45$ it seems that there are no more periodic orbits (and the rotation number is irrational). Thus here the small change of the parameter  by $0.05$ causes the real qualitative change in the behavior of the system. However, we must recall that usually (i.e. unless  the firing phase map is conjugated to the rational rotation), having the rational rotation number of a particular value $p/q$ is stable with respect to the small change of parameters and thus the system (in terms of periodic orbits) behaves in the same way,  which is what we call \emph{phase locking}. In fact the smaller the denominator of the rotation number, the more stable it is. Thus the small change of parameters within the neighborhood of the firing map with rational rotation number, also does not cause a qualitative change of the behaviour of the system, as long as this change of parameters preserves the rotation number (recall that, under some constraints, the mapping $\Phi\mapsto\varrho$ is a Devil-staircase, strictly increasing at irrational values and constant at rational ones with the longest intervals of being constant occurring at \emph{Farey fractions}, cf. e.g. Proposition 11.1.11 in \cite{katok}).

In Figure \ref{fig:beta2}  we may see also what happens for values of $\beta$ greater than $0.5$, precisely for $\beta=1$ and $\beta=2$. However,  for these parameter values  the firing map is not a homeomorphism any more and thus in particular the results might depend on the initial condition $(t,0)$. Therefore these cases are beyond the scope of this work.

\begin{remark} Note that the analogous statements of Propositions \ref{srednia}, \ref{slabazb} and \ref{przyblizaniedlalif} hold for the distribution of the firing phases, which is simply the invariant measure $\mu$.
\end{remark}
\section{Discussion}
We have shown many specific properties of the interspike-interval sequence arising from linear periodically driven integrate-and-fire models for which the emerging firing phase map is a circle homeomorphism. Many of them hold also for the general class $\dot{x}=F(t,x)$ of integrate-and-fire models (thus also non-linear models), as we prove in forthcoming paper ``Firing map and interspike-intervals for the general class of integrate-and-fire models with periodic drive''.

However, it would be interesting to have such rigorous results on interspike-intervals for periodically driven integrate-and-fire models with the firing phase map being not necessary a homeomorphism, but for instance just a continuous circle map. We predict that in such systems  greater variety of phenomena may be observed, mainly due to the fact that in this case we have \emph{rotation intervals} instead of the unique rotation number.

The natural extension of this research is detailed description of the interspike-interval sequence for IF systems with an almost periodic input (\cite{wmjs1}) and for bidimensional IF models (\cite{br-tob}).
\subsection*{Acknowledgements} The first author was supported  by national research grant NCN 2011/03/B\\/ST1/04533 and the second author by National Science Centre grant DEC-2011/01/N/ST1/02698.


%
\end{document}